\documentclass[11pt,oneside]{article}

\usepackage{amsthm,amsmath,amssymb,amsfonts}
\usepackage{mathtools,xfrac}
\usepackage{url,setspace,color}
\usepackage{subfigure,graphicx,epsfig,tikz,caption}
\usepackage{enumitem}
\usepackage{pdfsync,array}
\usepackage[normalem]{ulem}
\usepackage{algorithm,algorithmic}
\newcommand{\LINEIF}[2]{%
    \STATE\algorithmicif\ {#1}\ \algorithmicthen\ {#2}%
}

\usepackage{geometry,verbatim}
\usepackage{makeidx,latexsym}
\usepackage{bm}
\usepackage{booktabs}
\usepackage{multirow}
\usepackage{subfigure}

\usepackage[colorlinks=true]{hyperref}

\usepackage{authblk}
\usepackage[in]{fullpage}
\makeatletter
\def\imod#1{\allowbreak\mkern10mu({\operator@font mod}\,\,#1)}
\makeatother
\usepackage[compact]{titlesec}

\makeatletter
\renewenvironment{proof}[2][\proofname]{\par
\pushQED{\qed}
\normalfont
\topsep6\p@\@plus6\p@\relax
\trivlist
\item[\hskip\labelsep\textit{#2}\@addpunct{.}]\ignorespaces
}{%
\popQED\endtrivlist\@endpefalse
}
\makeatother

\def\Halmos{\mbox{\quad$\square$}}
\def\Halmos{\qed}
\def\Halmos{}

\theoremstyle{plain}
\newtheorem{theorem}{Theorem}[section]

\newtheorem{lemma}{Lemma}[section]
\newtheorem*{lemma*}{Lemma}
\newtheorem{proposition}{Proposition}[section]
\newtheorem*{proposition*}{Proposition}
\newtheorem{corollary}{Corollary}[section]

\newtheorem*{observation*}{Observation}
\newtheorem{assumption}{Assumption}[section]

\theoremstyle{definition}

\theoremstyle{remark}
\newtheorem{remark}{Remark}[section]
\newtheorem*{remark*}{Remark}

\newtheorem*{example*}{Example}

\def\TABLE{\caption}

\definecolor{gold}{rgb}{0.85,0.65,0}
\usepackage{xcolor}

\newcommand{\be}{\begin{eqnarray}}
\newcommand{\ee}[1]{\label{#1}\end{eqnarray}}
\newcommand{\ese}{\end{eqnarray*}}
\newcommand{\bse}{\begin{eqnarray*}}
\def\beq{\begin{equation}}
\def\eeq{\end{equation}}

\def\fnote#1{\footnote}
\newcommand{\epr}{\hfill\hbox{\hskip 4pt \vrule width 5pt height 6pt depth 1.5pt}\vspace{0.0cm}\par}
\def\qed{\ \hfill$\square$\par\smallskip}

\def\ra{\rangle}
\def\la{\langle}

\newcommand{\grad}{\ensuremath{\nabla}}

\def\E{{\mathbb{E}}}

\def\N{{\mathbb{N}}}

\def\R{{\mathbb{R}}}

\def\Se{{\mathbb{S}}}

\def\cA{{\cal A}}

\def\cM{{\cal M}}

\def\cP{{\cal P}}

\def\cS{{\cal S}}

\def\cV{{\cal V}}

\newcommand{\bbE}{\mathbb{E}}

\newcommand{\bbR}{\mathbb{R}}

\newcommand{\CA}{\mathcal{A}}

\newcommand{\CL}{\mathcal{L}}
\newcommand{\CM}{\mathcal{M}}

\newcommand{\CR}{\mathcal{R}}

\newcommand{\CV}{\mathcal{V}}

\DeclareMathOperator{\Opt}{Opt}

\DeclareMathOperator*{\argmin}{arg\,min}
\DeclareMathOperator*{\argmax}{arg\,max}

\DeclareMathOperator{\SadVal}{SV}

\DeclareMathOperator{\sad}{sad}

\DeclareMathOperator{\Prox}{Prox}

\def\epsilonsad{\epsilon_{\sad}}

\newcommand{\dgf}{d.g.f.}

\DeclareMathOperator{\sign}{sign}

\DeclareMathOperator{\Pess}{Pess}
\DeclareMathOperator{\Nom}{Nom}

\DeclareMathOperator{\Diag}{Diag}

\def\Fro{\mathop{\hbox{\rm\scriptsize Fro}}}
\def\Spec{\mathop{\hbox{\rm\scriptsize Spec}}}

\DeclareMathOperator{\cl}{cl}

\def\log{\mathop{{\rm log}}}

\def\square{\hbox{\vrule\vbox{\hrule\phantom{o}\hrule}\vrule}}

\begin{document}

\title{Online First-Order Framework for Robust Convex Optimization}
\author[1]{Nam Ho-Nguyen}
\author[1]{Fatma K{\i}l{\i}n\c{c}-Karzan}
\affil[1]{Tepper School of Business, Carnegie Mellon University, Pittsburgh, PA, 15213, USA.}
\date{July 20, 2016; revised October 31, 2016 and November 14, 2017}
\maketitle

\begin{abstract}

Robust optimization (RO) has emerged as one of the leading paradigms to efficiently model parameter uncertainty. The recent connections between RO and problems in statistics and machine learning domains demand for solving RO problems in ever more larger scale. However, the traditional approaches for solving RO formulations based on building and solving  robust counterparts or the iterative approaches utilizing nominal feasibility oracles can be prohibitively expensive and thus significantly hinder the scalability of RO paradigm. In this paper, we present a general and flexible iterative framework to approximately solve robust convex optimization problems that is built on a fully online first-order paradigm. In comparison to the existing literature, a key distinguishing feature of our approach is that it only requires access to first-order oracles that are remarkably cheaper than pessimization or nominal feasibility oracles, while maintaining the same convergence rates. This, in particular, makes our approach much more scalable and hence preferable in large-scale applications, specifically those from machine learning and statistics domains. 
We also provide new interpretations of existing iterative approaches in our framework and illustrate our framework on robust quadratic programming.

\end{abstract}

\maketitle

\section{Introduction}\label{sec:intro}

Robust optimization (RO) is one of the leading modeling paradigms for optimization problems under uncertainty. As opposed to the other approaches, RO seeks a solution that is immunized against \emph{all} possible realizations of uncertain model parameters (noises) from a given uncertainty set. It is widely adopted in practice mainly because of its computational tractability. We refer the reader to the paper by Ben-Tal and Nemirovski \cite{BenTalNem1998}, the book by Ben-Tal et al.\@ \cite{BenTalelGhaouiNemirovski2009} and surveys \cite{BenTalNemirovski2002,BenTalNemirovski2008,BertsimasBrownCaramanis2011,CaramanisMannorXu2011} for a detailed account of RO theory and numerous applications.

Recently, fascinating connections have been established between problems from the statistics and machine learning domains and robust optimization. More precisely, it is demonstrated that RO can be used to achieve desirable statistical properties such as stability, sparsity, and consistency. 
For example, for linear regression problems, El Ghaoui and Lebret \cite{ElGhaouiLebret1997} and Xu et al.\@ \cite{XuCaramanisMannor2010} respectively establish the equivalence of the ridge regression and Lasso to specific RO formulations of unregularized regression problems. Moreover, Xu et al.\@ \cite{XuCaramanisMannor2009} exhibit similar results in the context of regularizing support vector machines (SVMs), and \cite{XuCaramanisMannor2009,XuCaramanisMannor2010}  validate the statistical consistency of methods such as SVM and Lasso via RO methodology.

In addition to these RO interpretations of regularization techniques used in statistics and machine learning, robust versions of many problems from these domains are gaining traction. For example, \cite{Shivaswamy2006} examines robust variants of SVMs and other classification problems, and \cite{Ben-TalRobustSVM2012} explores a robust formulation for kernel classification problems. We refer the reader to \cite{CaramanisMannorXu2011,BenTalHazan2015} and references therein for further examples and details on connections between robust optimization and statistics and machine learning. 

These recent connections not only highlight the importance of RO methodology but also present algorithmic challenges where the scalability of RO algorithms with problem dimension becomes crucial. The primary method for solving a robust convex optimization problem is to transform it into an equivalent deterministic problem called the \emph{robust counterpart}. Under mild assumptions, this yields a convex and tractable robust counterpart problem  (see \cite{BenTalelGhaouiNemirovski2009,BertsimasBrownCaramanis2011,BenTalDenHertog2015}), which can then be solved using existing convex optimization software and tools. 
This traditional approach has seen much success in decision making domain, nevertheless it has a major drawback that the reformulated robust counterpart is often not as scalable as the deterministic nominal program. In particular, the robust counterpart can easily belong to a different class of optimization problems as opposed to the underlying original deterministic problem. For example, a linear program (LP) with ellipsoidal uncertainty is equivalent to a convex quadratic program (QP), and similarly, a conic-quadratic program with ellipsoidal uncertainty is equivalent to a semidefinite program (SDP) (see e.g., \cite{BenTalelGhaouiNemirovski2009, BertsimasBrownCaramanis2011}). It is well-known that convex QPs as opposed to LPs, and SDPs as opposed to convex QPs are much less scalable in practice. This then presents a critical challenge in applying RO methodology in big data applications frequently encountered in machine learning and statistics, where even solving the original deterministic nominal problem to high accuracy is prohibitively time-consuming. 

The iterative schemes that alternate between the generation/update of candidate solutions and the realizations of noises offer a convenient remedy to the scalability issues associated with the robust counterpart approach. Thus far, such approaches \cite{MutapcicBoyd2009} and \cite{BenTalHazan2015} have relied on two oracles: $(i)$ \emph{solution oracles} to solve  instances of extended (or nominal) problems with constraint structures similar  to (or the same as) the deterministic problem, and $(ii)$ \emph{noise oracles} to generate/update  particular realizations of the uncertain parameters. At each iteration of these schemes, both solution and noise oracles are called, and their outputs are used to update the inputs of each other oracle in the next iteration. Because solution oracles rely on a solver of the same class capable of solving the deterministic problem, these iterative approaches circumvent the issue of the robust counterpart approach potentially relying on a different solver. Nevertheless, these iterative approaches still suffer from a  serious drawback: the solution oracles in \cite{MutapcicBoyd2009,BenTalHazan2015} themselves can be expensive as they require solving  extended or nominal optimization problems completely. While solving the nominal problem is not as computationally demanding as solving the robust counterpart, the overall procedure depending on repeated calls to such oracles can be prohibitive. In fact, each such call to a solution oracle may endure a significant computational cost,  which is at least as much as the computational cost of solving an instance of the deterministic nominal problem. Note that, to ensure scalability, most applications in machine learning and statistics already need to rely on cheap first-order methods for solving deterministic nominal problems.

In this paper, we propose an efficient iterative framework for solving robust convex optimization problems which can rely on, in an \emph{online} fashion, much cheaper \emph{first-order oracles} in place of full solution and noise oracles. In particular, in each iteration, instead of solving a complete optimization problem within the solution and/or noise oracles, we show that simple simultaneous updates on the solution and noise in an online fashion using only first-order information from the deterministic constraint structure is sufficient to solve robust convex optimization problems. Moreover, we show that the number of calls to such online first-order (OFO)  oracles is not only at most that of the state-of-the-art iterative approaches utilizing full optimization based oracles for solution and/or noise, but also almost independent of the dimension of the problem. Therefore, this makes our approach especially attractive for applications in statistics and machine learning domains where it is critical to maintain that the overall approach has both gracious dependence on the dimension of the problem and cheap iterations. We outline our contribution more concretely after discussing the most relevant literature.

\subsubsection*{Related Work}

Thus far, the iterative approaches, which bypass the restrictions of the robust counterparts, work with \emph{extended} nominal problems that belong to the same class as the deterministic nominal one by carefully controlling the constraints included in the formulation corresponding to noise realizations. 

For robust binary linear optimization problems with only objective function uncertainty and a polyhedral uncertainty set, Bertsimas and Sim \cite{BertsimasSim2003} suggest an approach which relies on solving $n+1$ number of instances of the nominal problem, where $n$ is the dimension of the problem.

For robust convex optimization problems, Calafiore and Campi \cite{CalafioreCampi2004} study a `constraint sampling' approach based on forming a single extended nominal problem of the same class as the deterministic one via i.i.d.\@ sampling of noise realizations. They show that the optimal solution to this extended nominal problem is robust feasible with high probability where the probability depends on the sampling procedure, the number of samples drawn, and the dimension. 

Mutapcic and Boyd \cite{MutapcicBoyd2009} follow a `cutting-plane' type approach where in each iteration, a solution oracle is called to solve an extended nominal problem of the same class as the deterministic problem and a noise oracle, referred to as \emph{pessimization oracle}, is invoked to iteratively expand and refine the extended nominal problem. Given a candidate solution, a pessimization oracle either certifies its feasibility with respect to the robust constraints or returns a new noise realization from the uncertainty set for which the solution is infeasible; then the nominal constraint associated with that particular noise realization is included in the extended problem. This process is repeated until a robust feasible solution is found or the last extended problem is found to be infeasible. In the overall procedure, the number of iterations (or calls to the pessimization oracle) can be exponential in the dimension. Despite this, \cite{MutapcicBoyd2009} reports impressive computational results. The cutting-plane approach is also further tested on mixed integer linear problems in \cite{Bertsimas2016} and it is demonstrated that the same computational phenomenon holds.

Both of the approaches from \cite{CalafioreCampi2004} and \cite{MutapcicBoyd2009} pose issues for high-dimensional problems. 
In \cite{CalafioreCampi2004}, as the dimension grows, an extended problem with linearly more nominal constraints is required to ensure the high probability guarantee on finding a good quality solution. In \cite{MutapcicBoyd2009} at each iteration, a nominal constraint is added to the extended nominal problem. The theoretical bound on the number of constraints that need to be added is exponential, so the extended problem in  \cite{MutapcicBoyd2009} can grow to be exponentially large. Moreover, in both cases the extended nominal problem may no longer have  certain favorable problem structure of the deterministic nominal problem, such as a network flow structure.

To address these issues, in particular, the issue of solving extended nominal problems that are not only larger-in-size than the deterministic problem but also may lack certain favorable problem structure of the deterministic problem, Ben-Tal et al.\@ \cite{BenTalHazan2015} introduce a new iterative approach to approximately solve robust feasibility problems via a \emph{nominal feasibility oracle} and running an online learning algorithm to choose noise realizations. Given a particular noise realization, the nominal feasibility oracle solves an instance of  the deterministic nominal feasibility problem obtained by simply fixing the noise to the given value. Hence, the problem solved by this oracle has the \emph{same} number of constraints and the \emph{same structure} as the original nominal problem; in particular its size does not grow in each iteration. This is an important distinguishing feature of this approach. The other distinguishing feature is that Ben-Tal et al.\@ \cite{BenTalHazan2015} replace the pessimization oracle of \cite{MutapcicBoyd2009} by employing an online learning algorithm, which simply requires first-order information of the noise from the constraint functions. Moreover, \cite{BenTalHazan2015} provides a dimension independent bound on the number of iterations (nominal feasibility oracle calls).

Because the approaches of both \cite{MutapcicBoyd2009} and \cite{BenTalHazan2015} are closely related to our work, we give a detailed summary of these in Sections~\ref{sec:pessimization-oracle}~and~\ref{sec:oracle-based-approach} respectively and highlight their connections to our work. In fact,  we show that they both can be seen as special cases of our framework. 

We close with a brief summary of the assumptions on the computational requirements of these methods. The constraint sampling approach of \cite{CalafioreCampi2004} requires access to a sampling procedure on the uncertainty sets as well as an oracle capable of solving the extended nominal problem. The cutting plane approach of \cite{MutapcicBoyd2009} replaces the sampling procedure of \cite{CalafioreCampi2004} with a noise oracle, namely the pessimization oracle that works with the uncertainty sets but still requires the same type of optimization oracle as a solution oracle to solve the extended problems. Ben-Tal et al.\@ \cite{BenTalHazan2015} substitute the pessimization oracle with an online learning-based procedure, which requires merely first-order information from the constraint functions and simple projection type operations on the associated uncertainty sets, but it still relies on a solution oracle capable of solving the original nominal problem, which is essentially the same (up to log factors) as the optimization oracles in \cite{CalafioreCampi2004} and \cite{MutapcicBoyd2009}. If the deterministic problem admits special structure such as network flows etc., a specific solver can be used in the framework of \cite{BenTalHazan2015}, but this is not possible for \cite{CalafioreCampi2004} and \cite{MutapcicBoyd2009}.

\subsubsection*{Summary of Our Contributions}
It is possible to view all of these iterative approaches as two iterative processes that run simultaneously and in conjunction with each other to generate/update solutions and noise realizations. This naturally leads to a dynamic game environment where in each round Player 1 chooses a solution and Player 2 chooses a realization of uncertain parameters. In this framework, the policies employed by these players in their decision making determine the nature of the final approach. In the case of  \cite{MutapcicBoyd2009}, Player 1 considers all of the previous noise realizations when making his decision, whereas Player 2 simply reacts to the current solution when choosing the noise. In \cite{BenTalHazan2015}, Player 1 reacts to only the current noise in generating/updating the solution while Player 2 minimizes the regret associated with past solutions in choosing noise.

In this paper, we further analyze this interaction between Player 1 and Player 2, with the aim of deriving a simpler and computationally much less demanding iterative approach to solving RO problems. Our contributions can be summarized as follows.
\begin{enumerate}
\item We build a \emph{general and flexible framework} for iteratively solving robust feasibility problems, and demonstrate its flexibility by describing it as a meta-template. By customizing our framework appropriately, we modify the pessimization oracle-based approach of \cite{MutapcicBoyd2009} by replacing the extended nominal solver used in \cite{MutapcicBoyd2009} with efficient first-order updates. We call this the \emph{FO-based pessimization approach}, and demonstrate that as opposed to \cite{MutapcicBoyd2009} it  has both a much better bound on the number of oracle calls and far superior practical performance.  We also provide a new interpretation of the nominal feasibility oracle-based approach of \cite{BenTalHazan2015} as a special case within our framework. Furthermore, we extend the analysis of the approach of \cite{BenTalHazan2015} under the same assumptions, e.g., access to a nominal optimization oracle, and show that it can solve the robust optimization problem directly without relying on a binary search (see Remark \ref{rem:oracle-based-optimal-sol}).

\item When the original deterministic problem admits first-order oracles capable of providing gradient/subgradient information on each constraint function, we demonstrate that \emph{online first-order} (OFO) algorithms can be used to iteratively generate/update solutions and noise realizations simultaneously in an online manner leading to robust feasibility/infeasibility certificates within our framework. In contrast to the approaches of \cite{MutapcicBoyd2009} and \cite{BenTalHazan2015}, which rely on full nominal feasibility oracles to generate points, our OFO-based approach only requires simple update rules in each iteration and thus has much lower per-iteration cost. Besides, our noise oracle generates a realization of the noise in an online learning fashion as was done in \cite{BenTalHazan2015}, and hence it is less expensive than the pessimization oracle of \cite{MutapcicBoyd2009}. 

\item In our framework, the number of iterations (or oracle calls) needed to obtain approximate robust solution or a robust infeasibility certificate is a function of the approximation guarantee $\epsilon$ and the complexities of the domains for the solution and the uncertainty set; in particular, our convergence rate is (almost) independent of both the number of robust constraints and the dimension of the deterministic problem. We also demonstrate that the iteration complexity of our OFO-based approach is at least as good as that of the efficient approach of \cite{BenTalHazan2015}, and better than the exponential complexity of \cite{MutapcicBoyd2009}. Overall, our OFO-based approach leads to computational savings over the approach of \cite{BenTalHazan2015} by a factor as large as $O(1/(\epsilon^2 \log(1/\epsilon)))$ arithmetic operations when the number updates of the solution is smaller than or equal to the number of updates of the noise realization, which is the case in many applications. For further comparisons and discussion, see Section~\ref{sec:RateDiscussion}. In addition, our framework is amenable to exploiting favorable structural properties of the constraint functions such as strong concavity, smoothness, etc., through which better convergence rates can be achieved.

\item Our framework is based on formulating the robust feasibility problem as a convex-nonconcave saddle point (SP) problem, and explicitly analyzing its structure. While convex-concave SP problems are well-studied in the literature, and many efficient first-order algorithms exist for these (see for example \cite{Nesterov2005,JuditNem2012Pt1,JuditNem2012Pt2}), the convex-\emph{nonconcave} SP problem is not as well-studied. To our knowledge, an explicit study of convex-nonconcave SP problems and their relation to RO has not been conducted previously; in this respect, the most closely related work \cite{BenTalHazan2015} neither provides an explicit connection between robust feasibility and SP problems, nor analyzes their structure explicitly. 
\end{enumerate}

To demonstrate the application and effectiveness of our proposed framework, we walk through a detailed example on robust QPs. In particular, for robust QPs, we are able to leverage a recent convex QP-based reformulation of the classical trust region subproblem \cite{JeyakumarLi2013,Ho-NguyenKK2017TRS}
in order to avoid working with a nonconvex reformulation in a lifted space as in \cite[Section 4.2]{BenTalHazan2015} and relying on a probabilistic follow-the-perturbed-leader type algorithm \cite[Section 3.2]{BenTalHazan2015}. While using such nonconvex techniques will work within our framework, our convex reformulation allows us to work directly in the original space of the variables with a deterministic subgradient-based algorithm while still achieving asymptotically similar iteration complexity guarantees as \cite{BenTalHazan2015}. Moreover, each iteration of our approach requires only first-order updates where the most expensive operation is the computation of a maximum eigenvector; thus our per-iteration cost is significantly less. 

We also conduct a preliminary numerical study on the comparison of our approach with other iterative approaches \cite{MutapcicBoyd2009} and \cite{BenTalHazan2015} on robust QPs arising in portfolio optimization. Our results show that when the problem size is small, the nominal solver approaches of \cite{MutapcicBoyd2009} and \cite{BenTalHazan2015} are more efficient. However, when problem size increases, replacing the nominal solvers with first-order updates using our framework allows us to achieve faster solution times. This highlights the benefits and potential of investigating first-order based approaches such as ours in iterative RO methods.

\subsubsection*{Outline}
The rest of the paper is organized as follows. We begin with some notation and preliminaries in Section~\ref{sec:prelim}. We introduce our robust feasibility problem and robust feasibility/infeasibility certificates in Section~\ref{sec:R-feas-prelim}, convex-concave SP problems in Section~\ref{sec:SPprelim}, and briefly summarize important online convex optimization (OCO) tools as well as a useful OFO algorithm in Section~\ref{sec:OCO}. We formulate the robust feasibility problem as a convex-nonconcave SP problem in Section~\ref{sec:SPFormulation}; this formulation and certain bounds associated with its SP gap function form the basis of our general framework for solving robust feasibility problems. In Section~\ref{sec:OCOforRO} we specify an assortment of approaches obtained in our general framework by using different oracles. We examine our OFO-based approach in Section~\ref{sec:RegretforRO} by interpreting various terms in our framework in the context of OCO. In Section \ref{sec:pessimization-oracle}, we modify the pessimization oracle-based approach of \cite{MutapcicBoyd2009} to obtain an efficient bound on the number of iterations required. In Section \ref{sec:oracle-based-approach} we show how the nominal feasibility oracle-based approach of \cite{BenTalHazan2015} fits within our framework. 
Finally, we discuss the convergence rates and accelerations attainable in our framework and compare our work with the existing approaches in Section~\ref{sec:RateDiscussion}. In Section~\ref{sec:Application} we illustrate our OFO-based approach through an  example application on robust QPs. 
We provide in Section~\ref{sec:numerical} a preliminary numerical study comparing our framework with other iterative approaches \cite{MutapcicBoyd2009} and \cite{BenTalHazan2015}. We close with a summary of our results and a few compelling further research directions in Section~\ref{sec:Conclusions}. 
In Appendix~\ref{sec:concaveSP} we give an alternative formulation of the robust feasibility problem as a convex-concave SP problem in an extended space, and discuss its advantages and disadvantages over the convex-nonconcave SP formulation. 

\section{Notation and Preliminaries}\label{sec:prelim}

Given $a\in\R$, $\sign(a)$ denotes the sign of the number $a$. For a positive integer $n\in\N$,  we let $[n]=\{1,\ldots,n\}$ and define $\Delta_n:=\{x\in\R^n_+:~\sum_{i\in[n]} x_i=1\}$ to be the standard simplex. Throughout the paper, the superscript, e.g., $f^i, u^i, U^i$, is used to attribute items to the $i$-th constraint, whereas the subscript, e.g., $x_t,f_t,\phi_t$, is used to attribute items to the $t$-th iteration. Therefore, we sometimes use $u^i$, $x_t$, as well as $u^i_t$ to denote vectors in $\R^n$. We use the notation $\{x_t\}_{t=1}^T$ to denote the collection of items $\{x_1,\ldots,x_T\}$. Given a vector $x\in\R^n$, we let $x^{(k)}$ denote its $k$-th coordinate for $k\in[n]$. One exception we make to this notation is that we always  denote the convex combination weights $\theta\in\Delta_T$ with $\theta_t$. For $x\in\R^n$ and $p\in[1,\infty]$, we use $\|x\|_p$ to denote the $\ell_p$-norm of $x$ defined as  
\[
\|x\|_p=\begin{cases} \left(\sum_{i\in[n]} |x^{(i)}|^p \right)^{1/p} &\text{if } p\in[1,\infty) \\ \max_{i\in[n]}|x^{(i)}| &\text{if } p=\infty\end{cases} .
\]
Throughout this paper, we use Matlab notation to denote vectors and matrices, i.e., $[x;y]$ denotes the concatenation of two column vectors $x$, $y$. 
$\Se^n$ denotes the space of $n\times n$ symmetric matrices, and we let $\Se^n_+$ be the positive semidefinite cone in $\Se^n$. We let $I_n$ denote the identity matrix in $\Se^n$. For a matrix $A\in\Se^n$,  $\lambda_{\max}(A)$, $\|A\|_{\Fro}$, and $\|A\|_{\Spec}$ correspond to its maximal eigenvalue, Frobenius norm, and  spectral norm, respectively. 
Given a set $V$, we denote its closure by $\cl(V)$. We abuse notation slightly by denoting $\grad f(x)$ for both the gradient of function $f$ at $x$ if $f$ is differentiable and a subgradient of $f$ at $x$, even if $f$ is not differentiable. If $f$ is of the form $f(x,u)$, then $\grad_x f(x,u)$ denotes the subgradient of $f$ at $x$ while keeping the other variables fixed at $u$.

\subsection{Robust Feasibility Problem}\label{sec:R-feas-prelim}
Consider a convex \emph{deterministic} or \emph{nominal} mathematical program 
\begin{equation}\label{eqn:deterministic-opt}
\min_x \left\{ f^0(x):~ x \in X, \;~f^i(x,u^i) \leq 0,\ \forall i\in[m] \right\},
\end{equation}
where the domain $X\subset\R^n$ is closed and convex, the functions $f^0(x)$ and $f^i(x,u^i)$ for $i\in[m]$ are convex functions of $x$, and $u = (u^{1},\ldots,u^{m})$ is a fixed parameter vector. Without loss of generality we assume the objective function $f^0(x)$ does not have uncertainty. The \emph{robust convex optimization problem} associated with \eqref{eqn:deterministic-opt} is
\begin{equation}\label{eqn:robust-opt-intro}
\Opt:=\min_x \left\{ f^0(x):~ x \in X, \;~\sup_{u^i \in U^i} f^i(x,u^i) \leq 0,\ \forall i\in[m] \right\},
\end{equation}
where $U^1,\ldots,U^m$ are the \emph{uncertainty sets} given for the parameter $u^i$ of constraint $i\in[m]$. Because we assume formulation~\eqref{eqn:deterministic-opt} is convex, the overall optimization problem in \eqref{eqn:robust-opt-intro} is convex.

In this paper, we work under the following mild regularity assumption:
\begin{assumption}\label{ass:f-concave-u}
The constraint functions $f^i(x,u^i)$ for all $i\in[m]$ are finite-valued on the domain $X\times U^i$, convex in $x$ and concave in $u^i$. $X$, the domain for $x$, is closed and convex, and $U^i$, the domains for $u^i$, are closed and bounded.
\end{assumption}
We take Assumption~\ref{ass:f-concave-u} as given for all our results and proofs. Without loss of generality, we assume that the uncertainty set has a Cartesian product form $U^1 \times \ldots \times U^m$, see e.g., \cite{BenTalNemirovski2002}; we let $U = U^1 \times \ldots \times U^m$ and write $u = [u^1;\ldots;u^m] \in U$. We do not further assume that the sets $U^i$ are convex. However, for some algorithms we consider, convexity of $U^i$ for $i\in[m]$ will be required.

A convex optimization problem can be solved by solving a polynomial number of associated feasibility problems in a standard way, via a binary search over its optimal value. In particular, let $[\underline{\upsilon}_0,\overline{\upsilon}_0]$ be an initial interval containing the optimal value of \eqref{eqn:robust-opt-intro}. At each iteration $k$ of the binary search, we update the domain $X_k:=X\cap\{x:\;f^0(x)\leq \upsilon_k$\} for some $\upsilon_k\in[\underline{\upsilon}_k,\overline{\upsilon}_k]$ and arrive at the following robust feasibility problem:
\begin{equation}\label{eqn:robust-feas}
\text{find}\ \ x \in X_k \quad\text{s.t.}\quad \sup_{u^i\in U^i} f^i(x,u^i) \leq 0 \quad \forall i\in[m].
\end{equation}
Then based on the feasibility/infeasibility status of \eqref{eqn:robust-feas}, we update our range $[\underline{\upsilon}_{k+1},\overline{\upsilon}_{k+1}]$ and go to iteration $k+1$. In this scheme, we are guaranteed to find a solution $x^* \in X$ whose objective value is within $\delta>0$ of the optimum value of \eqref{eqn:robust-opt-intro} in at most $\left\lfloor \log_2\left(\frac{\overline{\upsilon}_0-\underline{\upsilon}_0}{\delta}\right) \right\rfloor$ iterations. Therefore, one can equivalently study the complexity of solving  robust feasibility problem~\eqref{eqn:robust-feas} as opposed to \eqref{eqn:robust-opt-intro}. 
From now on, we focus on solving robust feasibility problem and assume that the constraint on the objective function $f^0(x)$ is already included in the domain $X$ for simplicity in our notation.

Given functional constraints $f^i(x) \leq 0$, $i \in [m]$, most convex optimization methods will declare infeasibility or return an approximate solution $x \in X$ such that $f^i(x) \leq \epsilon$ for $i \in [m]$ for some tolerance level $\epsilon > 0$. Therefore, we consider the following \emph{robust approximate feasibility problem}:
\begin{equation}\label{eqn:approx-robust-feas}
\begin{cases}
\text{\emph{Either}: find}\ \ x \in X \quad\text{s.t.}\quad \sup_{u^i\in U^i} f^i(x,u^i) \leq \epsilon \quad \forall i\in[m];\\
\text{\emph{or}: declare infeasibility, } \forall x \in X,\ \exists i \in [m] \quad \text{s.t.} \quad \sup_{u^i\in U^i} f^i(x,u^i) > 0.
\end{cases} 
\end{equation}

We refer to any feasible solution $x$ to \eqref{eqn:approx-robust-feas}, i.e., $x \in X$ such that $\sup_{u^i\in U^i} f^i(x,u^i) \leq \epsilon$ holds for all $i\in[m]$ as a \emph{robust $\epsilon$-feasibility certificate}. Similarly, any realization of the uncertain parameters $\bar{u}\in U$ such that there exists no $x\in X$ satisfying $f^i(x,\bar{u}^i) \leq 0$  for all $i\in[m]$ is referred to as a \emph{robust infeasibility certificate}.

\subsection{Saddle Point Problems}\label{sec:SPprelim}

Saddle point (SP) problems play a vital role in our developments. In its most general form, a convex-concave SP problem is given by
\begin{equation}\label{eqn:SadVal}
\SadVal=\inf_{x \in X} \sup_{y \in Y} \phi(x,y),  \tag{$\cS$}
\end{equation}
where the function $\phi(x,y)$ is convex in $x$ and concave in $y$ and the domains $X,Y$ are nonempty closed convex sets in Euclidean spaces $\E_x,\E_y$. 

Any convex-concave SP problem \eqref{eqn:SadVal} gives rise to two convex optimization problems that are dual to each other:
\begin{equation*}\label{neq1}
\begin{array}{rclcr}
\Opt(P)&=&\inf_{x\in X}[ \overline{\phi}(x):=\sup_{y\in Y} \phi(x,y)]&&(P)\\
\Opt(D)&=&\sup_{y\in Y}[ \underline{\phi}(y):=\inf_{x\in X} \phi(x,y)] &&(D)\\
\end{array}
\end{equation*}
with $\Opt(P)=\Opt(D)=\SadVal$.
It is well-known that the solutions to \eqref{eqn:SadVal} --- the saddle points of $\phi$ on $X\times Y$ --- are exactly the pairs $[x;y]$ formed by optimal solutions to the problems $(P)$ and $(D)$.

We quantify the accuracy of a candidate solution $[\bar{x},\bar{y}]$ to SP problem \eqref{eqn:SadVal} with the \emph{saddle point gap} given by 
\begin{equation}\label{eqn:SPgap}
\epsilonsad^\phi(\bar{x},\bar{y}):=\overline{\phi}(\bar{x})-\underline{\phi}(\bar{y}) 
=\underbrace{\left[\overline{\phi}(\bar{x})-\Opt(P)\right]}_{\geq0}+
\underbrace{\left[\Opt(D)-\underline{\phi}(\bar{y})\right]}_{\geq0}.
\end{equation}

Because convex-concave SP problems are simply convex optimization problems, they can in principle be solved by polynomial-time interior point methods (IPMs).  However, the computational complexity of such methods depends heavily on the dimension of the problem. Thus, scalability of resulting algorithms becomes an issue in large-scale applications. 
As a result, for large-scale SP problems, one has to resort to first-order subgradient-type methods.
On a positive note, there are many efficient first-order methods (FOMs) for convex-concave SP problems. These in particular include Nesterov's accelerated gradient descent algorithm \cite{Nesterov2005} and Nemirovski's Mirror-Prox algorithm \cite{Nemirovski2005}, both of which bound the saddle point gap at a rate of $\epsilonsad^\phi(\bar{x}_T,\bar{u}_T) \leq O\left({1\over T}\right)$ where $\bar{x}_T,\bar{u}_T$ are solutions obtained after $T$ iterations.

\subsection{Online Convex Optimization Tools}\label{sec:OCO}

Our efficient framework for RO employs tools from the online convex optimization domain. We now briefly outline these and refer to \cite{CesaBianchiLugosi2006,Hazan2011,Shalev-Shwarz2011} for further details and applications of OCO.

OCO is used to capture decision making in dynamic environments. We are given a finite time horizon $T$, closed, bounded, and convex domain $Z$, and in each time period $t \in [T]$, a convex loss function $f_t:Z \to \R$ is revealed. At time periods $t \in [T]$ we must choose a decision $z_t \in Z$, and based on this we suffer a loss of $f_t(z_t)$ and receive some feedback typically in the form of first-order information on $f_t$. Our goal is to minimize the \emph{weighted regret}
\begin{equation}\label{eqn:weighted-regret-defn}
\sum_{t=1}^{T} \theta_t f_t(z_t) - \inf_{z\in Z} \sum_{t=1}^{T} \theta_t f_t(z),
\end{equation}
where $\theta \in \Delta_T$ is a vector of convex combination weights.\footnote{Note that in the OCO literature, regret is usually defined with uniform weights $\theta_t = 1/T$. Nonuniform weights introduce flexibility to our framework by allowing selection of specific customization of OCO algorithms for exploiting structural properties of the constraint functions $f^i$ to achieve better convergence rates. A prime example for this is when the functions are strongly convex.}

Most OCO algorithms are closely related to offline iterative FOMs. In this paper, we will make use of the proximal setup of \cite{JuditNem2012Pt1} to choose the sequence $\{z_t\}_{t=1}^T$ which ensures that the weighted regret \eqref{eqn:weighted-regret-defn} converges to $0$ as $T \to \infty$. Thus, we make the following assumption on $Z$ for the existence of a proximal setup.  
\begin{assumption}\label{ass:proximal-setup-OFO}
Let $\bbE_z$  be the Euclidean space containing $Z$. There exists a norm $\|\cdot\|$ and its dual norm $\|\cdot\|_*$ on $\bbE_z$, a distance-generating function $\omega:Z \to \R$ which is $1$-strongly convex with respect to $\|\cdot\|$ and leads to an easy-to-compute prox function $\Prox_z(\xi):=\argmin_{w\in Z} \left\{\la \xi,w\ra+\omega(w)-\la\omega'(z),w-z\ra\right\}$ and set width $\Omega := \max_{z \in Z} \omega(z) - \min_{z \in Z} \omega(z)$ which is finite when $Z$ is bounded.
\end{assumption}
The proximal setup of Assumption~\ref{ass:proximal-setup-OFO} allows us to adjust to the geometry of domain $Z$. The standard basic domains satisfying Assumption~\ref{ass:proximal-setup-OFO} include simplex, Euclidean ball, and spectahedron; see \cite[Section 1.7]{JuditNem2012Pt1} for the standard proximal setups (i.e., Assumption~\ref{ass:proximal-setup-OFO}) for these basic domains in terms of selection of $\|\cdot\|$ and resulting $\omega$, $\Prox$ computation, and set width $\Omega$. 

Under Assumption~\ref{ass:proximal-setup-OFO} and various structural properties, the straightforward extension of the standard online mirror descent algorithm (see, e.g., \cite{KakadeSSTewari2012}) from uniform weigths to weighted regret achieves the following convergence rate.

\begin{theorem}[{\cite[Theorem 5]{KakadeSSTewari2012}}]\label{thm:OCO-non-smooth}
Suppose there exists $G \in (0,\infty)$ such that $\|\grad f_t(z)\|_* \leq G$ for all $z \in Z$, $t \in [T]$. Define $\gamma = \sqrt{\frac{2 \Omega}{G^2 \sum_{t=1}^T \theta_t^2}}$. Choose $z_1 = \argmin_{z \in Z} \omega(z)$ and $z_{t+1} = \Prox_{z_t} (\gamma \theta_t \grad f_t(z_t))$ for $t\in[T]$. Then
\[ \sum_{t=1}^{T} \theta_t f_t(z_t) - \inf_{z\in Z} \sum_{t=1}^{T} \theta_t f_t(z) \leq \sqrt{2 \Omega G^2 \sum_{t=1}^T \theta_t^2}. \]
In particular, for uniform weights $\theta_t = 1/T$, the upper-bound becomes $O(1/\sqrt{T})$.
\end{theorem}
We refer to \cite{KakadeSSTewari2012} for details of the proof. When $\omega(z) = z^\top z/2$ and weights $\theta_t=1/T$ for $t\in[T]$, the update rule $z_{t+1} = \Prox_{z_t} (\gamma \grad f_t(z_t))$ becomes simply gradient descent, and Theorem~\ref{thm:OCO-non-smooth} reduces to the standard bound of online gradient descent from \cite{Zinkevich2003}.

\section{General Framework for Robust Feasibility Problems}\label{sec:SPFormulation}
In this section, we build a general framework to solve the robust feasibility problem \eqref{eqn:approx-robust-feas} by working with its natural saddle point formulation.

Given constraint functions $f^i(x,u^i)$, $i \in [m]$, let us define $\Phi(x,u):= \max_{i\in[m]} f^i(x,u^i)$. Then $\Phi(x,u)$ is a convex function of $x$, but not necessarily concave in $u$. In addition, with this definition of $\Phi(\cdot)$, the robust approximate feasibility problem~\eqref{eqn:approx-robust-feas} is equivalent to simply verifying  
\begin{equation}\label{eqn:convex-nonconcaveSP}
\text{either}\quad \inf_{x \in X} \sup_{u \in U} \Phi(x,u) = \inf_{x \in X} \max_{i\in[m]} \sup_{u^i \in U^i} f^i(x,u^i) \leq \epsilon \quad \text{or} \quad \inf_{x \in X} \sup_{u \in U} \Phi(x,u) > 0,
\end{equation}
which is nothing but solving a specific SP problem and checking its value. Analogous to the convex-concave SP gap \eqref{eqn:SPgap}, for a  given solution $[\bar{x},\bar{u}]$, we define the SP gap of problem \eqref{eqn:convex-nonconcaveSP} as 
\[
\epsilonsad^\Phi(\bar{x},\bar{u}):=\overline{\Phi}(\bar{x})-\underline{\Phi}(\bar{u}) 
=\sup_{u\in U} \Phi(\bar{x},u) - \inf_{x\in X} \Phi(x,\bar{u}).
\]

In general, solving a convex-nonconcave SP problem of form~\eqref{eqn:convex-nonconcaveSP}, i.e., finding a solution $[\bar{x},\bar{u}]$ such that $\epsilonsad^\Phi(\bar{x},\bar{u})\leq\epsilon$, can be difficult. That said, a bound on the SP gap $\epsilonsad^\Phi(\bar{x},\bar{u})$ along with the value of $\Phi(\bar{x},\bar{u})$ leads to robust  feasibility certificates for \eqref{eqn:convex-nonconcaveSP} as follows.
\begin{theorem}\label{thm:saddlept-certificate}
Let $\Psi:X \times U \to \R$ be a given function associated with a SP (not necessarily admitting a convex-concave structure). Suppose we have $\bar{x} \in X$, $\bar{u} \in U$, and $\tau\in(0,1)$ such that $\epsilonsad^\Psi(\bar{x},\bar{u})\leq \tau\epsilon$. 
Then if $\Psi(\bar{x},\bar{u}) \leq (1-\tau)\epsilon$, we have  $\sup_{u \in U} \Psi(\bar{x},u) \leq \epsilon$. Moreover, if $\Psi(\bar{x},\bar{u}) > (1-\tau)\epsilon$ and $\tau\leq{1\over 2}$, we have  $\inf_{x \in X} \Psi(x,\bar{u}) > 0$.
\end{theorem}
\proof{Proof.}
Suppose $\Psi(\bar{x},\bar{v}) \leq (1-\tau)\epsilon$. Because $\epsilonsad^\Psi(\bar{x},\bar{v}) =\sup_{u \in U} \Psi(\bar{x},u) - \inf_{x\in X} \Psi(x,\bar{u}) \leq \tau\epsilon$, we have $\sup_{u \in U} \Psi(\bar{x},u) \leq \inf_{x\in X} \Psi(x,\bar{u}) + \tau\epsilon \leq \Psi(\bar{x},\bar{u}) + \tau\epsilon \leq \epsilon$. On the other hand, when $\Psi(\bar{x},\bar{u}) > (1-\tau)\epsilon$, we have $(1-\tau)\epsilon < \Psi(\bar{x},\bar{u}) \leq \sup_{u \in U} \Psi(\bar{x},u) \leq \inf_{x \in X} \Psi(x,\bar{u}) + \tau\epsilon$, which implies $\inf_{x\in X} \sup_{u\in U} \Psi(x,u) \geq \inf_{x\in X} \Psi(x,\bar{u}) > (1-2\tau)\epsilon\geq 0$ when $\tau\leq{1\over 2}$.
\Halmos
\endproof

\begin{remark}\label{rem:SPm=1}
When $m=1$, $\Phi(x,u)=f^1(x,u^1)$, and it is thus convex in $x$ and concave in $u$ due to Assumption~\ref{ass:f-concave-u}. Therefore, in the case of a single robust constraint, i.e., $m=1$, under Assumption~\ref{ass:f-concave-u} and assuming $U=U^1$ is a closed convex set, the optimization problem in \eqref{eqn:convex-nonconcaveSP} reduces to a standard convex-concave SP problem. 
\epr
\end{remark}

While it is not very common, a few robust convex optimization problems come with a single robust constraint and convex uncertainty set $U$; see for example \cite{Ben-TalRobustSVM2012} for a robust version of a SVM problem with one constraint. In such cases, based on Remark~\ref{rem:SPm=1}, the resulting convex-concave SP problems can directly be solved via efficient FOMs. On the other hand, in the presence of multiple constraints, the function $\Phi(x,u)$ is not concave in $u = [u^1;\ldots;u^m]$ even under Assumption~\ref{ass:f-concave-u}. Nevertheless, when $m>1$, it is still possible to have a convex-concave SP reformulation of the optimization problem in \eqref{eqn:convex-nonconcaveSP} in an extended space via perspective transformations, which we present in Appendix~\ref{sec:concaveSP}. While this reformulation has the benefit of reducing the robust feasibility problem to a well-known and well-studied problem, it destroys the simplicity of the original domains and constraint functions and hence comes with some challenges. Therefore, we develop a framework where we work directly with the convex-nonconcave SP formulation in \eqref{eqn:convex-nonconcaveSP} in the space of original variables. 
Moreover, because we work in the original space of variables, we simply utilize the first-order information on the original constraint functions $f^i$ and original domains $X$ and $U^i$. This direct approach in particular allows us to take greater advantage of the structure of the original formulation such as the availability of efficient projection (prox) computations over domains $X,\, U^i$, and/or better parameters for smoothness, Lipschitz continuity, etc., of the functions $f^i$.

Because $\Phi(x,u)$ is not concave in $u$, we cannot bound the SP  gap $\epsilonsad^\Phi(\bar{x},\bar{u})$ by using traditional FOMs designed for solving convex-concave SP problems. However, we next show that by just \emph{partially} upper bounding $\epsilonsad^\Phi(\bar{x},\bar{u})$, we can derive a general iterative framework to obtain robust feasibility/infeasibility certificates. We describe the further specifics of this framework in Section~\ref{sec:OCOforRO}. 

Henceforth we will no longer use the shorthand notation $\Phi(x,u) = \max_{i \in [m]} f^i(x,u^i)$,  but we will denote the SP gap $\epsilonsad^\Phi(\bar{x},\bar{u})$ as
\begin{equation}\label{eqn:SVgap}
\epsilon(\bar{x},\bar{u}) := \epsilonsad^\Phi(\bar{x},\bar{u}) = \max_{i\in[m]} \sup_{u^i \in U^i} f^i(\bar{x},u^i) - \inf_{x\in X} \max_{i\in[m]} f^i(x,\bar{u}^i).
\end{equation}

The robust feasibility certificate result from Theorem~\ref{thm:saddlept-certificate} indicates the importance of bounding the SP gap $\epsilon(\bar{x},\bar{u})$. Often, FOMs achieve this by iteratively generating points $x_t \in X$, $u_t \in U$ for $t\in [T]$ and tracking the points $\bar{x}$ and $\bar{u}$ obtained from a convex combination of $\{x_t,u_t\}_{t=1}^T$. In order to simplify our notation, given convex combination weights $\theta\in\Delta_T$ and points $\{x_t,u_t\}_{t=1}^T$, we let 
\[ \bar{x}_T := \sum_{t=1}^{T} \theta_t x_t \quad\text{and}\quad \bar{u}_T := \sum_{t=1}^{T} \theta_t u_t. \]
We now present an upper bound on $\epsilon(\bar{x}_T,\bar{u}_T)$ that follows naturally from the convex-concave structure of functions $f^i$. 
To this end, given a set of vectors  $y_t\in\Delta_m$ for $t\in[T]$, we also define
\begin{align*}
\epsilon^\circ(\{x_t,u_t,\theta_t\}_{t=1}^T) &:= \max_{i\in[m]} \left\{\sup_{u^i\in U^i} \sum_{t=1}^{T} \theta_t f^i(x_t,u^i) - \sum_{t=1}^{T} \theta_t f^i(x_t,u_t^i)\right\},\quad\text{and}\\
\epsilon^\bullet(\{x_t,u_t,y_t,\theta_t\}_{t=1}^T) &:= \max_{i\in[m]} \sum_{t=1}^{T} \theta_t f^i(x_t,u_t^i) - \inf_{x\in X} \sum_{t=1}^{T} \theta_t \sum_{i=1}^{m} y^{(i)}_t f^i(x,u_t^i),
\end{align*}
together with 
\[
\widehat{\epsilon}(\{x_t,u_t,y_t,\theta_t\}_{t=1}^T) :=  \inf_{x\in X} \sum_{t=1}^{T} \theta_t \sum_{i=1}^{m} y^{(i)}_t f^i(x,u_t^i) - \inf_{x\in X} \max_{i\in[m]} \sum_{t=1}^{T} \theta_t f^i(x,u_t^i).
\]

Our next result relates these quantities to the value of the SP gap function $\epsilon\left(\bar{x}_T, \bar{u}_T\right)$.
\begin{proposition}\label{prop:gap-nonconcave}
Let $x_t \in X$ and $u_t \in U$ for $t\in[T]$ be given a set of vectors. Then for any set of vectors $y_t\in\Delta_m$ for $t\in[T]$ and any $\theta\in\Delta_T$, we have
\begin{equation}\label{eqn:saddlept-gap-bound}
\epsilon\left(\sum_{t=1}^{T} \theta_t x_t, \sum_{t=1}^{T} \theta_t u_t\right)  \leq \epsilon^\circ(\{x_t,u_t,\theta_t\}_{t=1}^T) +\epsilon^\bullet(\{x_t,u_t,y_t,\theta_t\}_{t=1}^T) + \widehat{\epsilon}(\{x_t,u_t,y_t,\theta_t\}_{t=1}^T).
\end{equation}
\end{proposition}

\proof{Proof.}
Given  $y_t\in\Delta_m$ for $t\in[T]$ and $\theta\in\Delta_T$, let us define $\bar{x} := \sum_{t=1}^{T} \theta_t x_t$ and $\bar{u} := \sum_{t=1}^{T} \theta_t u_t$. We first partition $\epsilon(\bar{x},\bar{u})$ as $\epsilon(\bar{x},\bar{u})=\overline{\epsilon}(\bar{x})+\underline{\epsilon}(\bar{u})$ where 
\begin{align*}
\overline{\epsilon}(\bar{x}) &:=\max_{i\in[m]} \sup_{u^i\in U^i} f^i(\bar{x},u^i) - \inf_{x\in X} \max_{i\in[m]} \sup_{u^i\in U^i} f^i(x,u^i) , \\
\underline{\epsilon}(\bar{u}) &:= \inf_{x\in X} \max_{i\in[m]} \sup_{u^i\in U^i} f^i(x,u^i) - \inf_{x\in X} \max_{i\in[m]} f^i(x,\bar{u}^i),
\end{align*}
and then derive upper bounds on $\overline{\epsilon}(\bar{x})$ and $\underline{\epsilon}(\bar{u})$. 

We start with bounding $\overline{\epsilon}(\bar{x})$. Because the functions $f^i(x,u^i)$ are convex in $x$ for all $i$ and $\theta\in\Delta_T$, we have $\max_{i\in[m]} \sup_{u^i\in U^i} f^i(\bar{x},u^i) \leq \max_{i\in[m]} \sup_{u^i\in U^i} \sum_{t=1}^{T} \theta_t f^i(x_t,u^i) $. Therefore, 
{\small
\begin{align}
\overline{\epsilon}(\bar{x}) &= \max_{i\in[m]} \sup_{u^i\in U^i} f^i(\bar{x},u^i)
- \inf_{x\in X} \max_{i\in[m]} \sup_{u^i\in U^i} f^i(x,u^i) \nonumber\\
 &\leq \max_{i\in[m]} \sup_{u^i\in U^i} \sum_{t=1}^{T} \theta_t f^i(x_t,u^i) - \max_{i\in[m]} \sum_{t=1}^{T} \theta_t f^i(x_t,u_t^i) + \max_{i\in[m]} \sum_{t=1}^{T} \theta_t f^i(x_t,u_t^i) - \inf_{x\in X} \max_{i\in[m]} \sup_{u^i\in U^i} f^i(x,u^i)\nonumber\\
 &\leq \max_{i\in[m]} \left\{\sup_{u^i\in U^i} \sum_{t=1}^{T} \theta_t f^i(x_t,u^i) - \sum_{t=1}^{T} \theta_t f^i(x_t,u_t^i)\right\} + \max_{i\in[m]} \sum_{t=1}^{T} \theta_t f^i(x_t,u_t^i) - \inf_{x\in X} \max_{i\in[m]} \sup_{u^i\in U^i} f^i(x,u^i) \label{eqn:epsOver},
\end{align}}
where the last inequality follows since $\max_{i\in[m]} \{\alpha_i-\beta_i\} \geq \max_{i\in[m]}\alpha_i - \max_{i\in[m]} \beta_i$ for any sequence of numbers $\alpha_i,\beta_i$, $i \in [m]$.

Note that $ \inf_{x\in X} \max_{i\in[m]} f^i(x,u^i) \geq  \inf_{x\in X} \max_{i\in[m]} \sum_{t=1}^{T} \theta_t f^i(x,u^i)$ because under Assumption~\ref{ass:f-concave-u} the functions $f^i(x,u^i)$ are concave in $u^i$ for all $i$.  
Thus, we arrive at
\begin{align}
\underline{\epsilon}(\bar{u}) &= \inf_{x\in X} \max_{i\in[m]} \sup_{u^i\in U^i} f^i(x,u^i) - \inf_{x\in X} \max_{i\in[m]} f^i(x,\bar{u}^i) \nonumber\\
 &\leq \inf_{x\in X} \max_{i\in[m]} \sup_{u^i\in U^i} f^i(x,u^i) - \inf_{x\in X} \sum_{t=1}^{T} \theta_t \sum_{i=1}^{m} y^{(i)}_t f^i(x,u_t^i) + \inf_{x\in X} \sum_{t=1}^{T} \theta_t \sum_{i=1}^{m} y^{(i)}_t f^i(x,u_t^i) \nonumber\\
 &\qquad\qquad - \inf_{x\in X} \max_{i\in[m]} \sum_{t=1}^{T} \theta_t f^i(x,u_t^i) \nonumber\\
 &= \inf_{x\in X} \max_{i\in[m]} \sup_{u^i\in U^i} f^i(x,u^i) - \inf_{x\in X} \sum_{t=1}^{T} \theta_t \sum_{i=1}^{m} y^{(i)}_t f^i(x,u_t^i) + \widehat{\epsilon}(\{x_t,u_t,y_t,\theta_t\}_{t=1}^T). \label{eqn:epsUnder}
\end{align}
Then by summing \eqref{eqn:epsOver} and \eqref{eqn:epsUnder} and rearranging the terms, we deduce the result.
\Halmos
\endproof

We are now ready to state our main result. This is analogous to Theorem~\ref{thm:saddlept-certificate} except that we do not need to bound all three terms in \eqref{eqn:saddlept-gap-bound}, but instead it suffices to guarantee that
\[ \epsilon^\circ(\{x_t,u_t,\theta_t\}_{t=1}^T) +\epsilon^\bullet(\{x_t,u_t,y_t,\theta_t\}_{t=1}^T) \leq \epsilon. \]
We show that when the above condition holds, based on the value of $\max_{i\in[m]} \sum_{t=1}^{T} \theta_t f^i(x_t,u_t^i)$ we can then obtain either a robust $\epsilon$-feasible solution, or an infeasibility certificate.

\begin{theorem}\label{thm:robust-feas-oracle}
Suppose we have sequences $\{x_t,u_t,y_t,\theta_t\}_{t=1}^T$ with $x_t\in X$, $u_t\in U$, $y_t\in\Delta_m$ for all $t\in[T]$, $\theta\in\Delta_T$. Let $\tau \in (0,1)$. 
If $\epsilon^\circ(\{x_t,u_t,\theta_t\}_{t=1}^T) \leq \tau \epsilon$ and $\max_{i\in[m]} \sum_{t=1}^{T} \theta_t f^i(x_t,u_t^i)\leq (1-\tau) \epsilon$, then the solution $\bar{x}_T := \sum_{t=1}^{T} \theta_t x_t$ is $\epsilon$-feasible with respect to \eqref{eqn:approx-robust-feas}. 
If $\epsilon^\bullet(\{x_t,u_t,y_t,\theta_t\}_{t=1}^T)\leq (1-\tau) \epsilon$ and $\max_{i\in[m]} \sum_{t=1}^{T} \theta_t f^i(x_t,u_t^i) > (1-\tau) \epsilon$, then \eqref{eqn:approx-robust-feas} is infeasible.
\end{theorem}
\proof{Proof.}
First suppose there exists a $\tau\in(0,1)$ and corresponding vectors $\{x_t,u_t,y_t,\theta_t\}_{t=1}^T$ such that $\epsilon^\circ(\{x_t,u_t,\theta_t\}_{t=1}^T) \leq \tau \epsilon$ and $\max_{i\in[m]} \sum_{t=1}^{T} \theta_t f^i(x_t,u_t^i) \leq (1-\tau) \epsilon$ holds as well. Note that
\begin{align}
\tau \epsilon &\geq \epsilon^\circ(\{x_t,u_t,\theta_t\}_{t=1}^T)
=\max_{i\in[m]} \left\{\sup_{u^i\in U^i} \sum_{t=1}^{T} \theta_t f^i(x_t,u^i) - \sum_{t=1}^{T} \theta_t f^i(x_t,u_t^i)\right\}\nonumber\\
&\geq \max_{i\in[m]} \sup_{u^i\in U^i} \sum_{t=1}^{T} \theta_t f^i(x_t,u^i) - \max_{i\in[m]} \sum_{t=1}^{T} \theta_t f^i(x_t,u_t^i) \label{eqn:thm:feas} ,
\end{align}
where the last inequality follows since $\max_{i\in[m]} \{\alpha_i-\beta_i\} \geq \max_{i\in[m]}\alpha_i - \max_{i\in[m]} \beta_i$ for any sequence of numbers $\alpha_i,\beta_i$, $i \in [m]$. Then $\bar{x}_T$ is an $\epsilon$-feasible solution for \eqref{eqn:approx-robust-feas} because 
\begin{align*}
\max_{i\in[m]} \sup_{u^i\in U^i} f^i\left(\bar{x}_T,u^i\right) &= \max_{i\in[m]} \sup_{u^i\in U^i} f^i\left(\sum_{t=1}^{T} \theta_t x_t,u^i\right) \\
&\leq \max_{i\in[m]} \sup_{u^i\in U^i} \sum_{t=1}^{T} \theta_t f^i(x_t,u^i) \leq \tau \epsilon + \max_{i\in[m]} \sum_{t=1}^{T} \theta_t f^i(x_t,u_t^i) \leq \epsilon,
\end{align*}
where the first inequality follows from the convexity of the functions $f^i$ and the fact that $\theta\in\Delta_T$, the second inequality from \eqref{eqn:thm:feas}, and the last inequality holds since $\max_{i\in[m]} \sum_{t=1}^{T} \theta_t f^i(x_t,u_t^i) \leq (1-\tau) \epsilon$. 

On the other hand, suppose $\epsilon^\bullet(\{x_t,u_t,y_t,\theta_t\}_{t=1}^T)\leq (1-\tau) \epsilon$ and $\max_{i\in[m]} \sum_{t=1}^{T} \theta_t f^i(x_t,u_t^i) > (1-\tau) \epsilon$. Note that
\begin{align}\label{eqn:epsPos}
\inf_{x\in X} \sum_{t=1}^{T} \theta_t \sum_{i=1}^{m} y^{(i)}_t f^i(x,u_t^i) &\leq \inf_{x\in X} \sum_{t=1}^{T} \theta_t \max_{i\in[m]} f^i(x,u_t^i) \nonumber\\
&\leq \inf_{x\in X} \sum_{t=1}^{T} \theta_t \max_{i\in[m]} \sup_{u^i\in U^i} f^i(x,u^i) 
= \inf_{x\in X} \max_{i\in[m]} \sup_{u^i\in U^i} f^i(x,u^i),
\end{align}
where the first inequality follows since $y_t\in\Delta_m$ for all $t\in[T]$, the second inequality holds because $f^i(x,u_t^i) \leq\sup_{u^i\in U^i} f^i(x,u^i)$ for all $i\in[m]$ and $y_t^{(i)} \geq 0$ for $i \in [m]$, $t \in [T]$, and the last equation follows from $\theta\in\Delta_T$. Then using the bound 
\begin{equation}\label{eqn:epsbullet-bound}
(1-\tau) \epsilon \geq \epsilon^\bullet(\{x_t,u_t,y_t,\theta_t\}_{t=1}^T) 
=\max_{i\in[m]} \sum_{t=1}^{T} \theta_t f^i(x_t,u_t^i) - \inf_{x\in X} \sum_{t=1}^{T} \theta_t \sum_{i=1}^{m} y^{(i)}_t f^i(x,u_t^i),
\end{equation}
we arrive at
\[
\inf_{x\in X} \max_{i\in[m]} \sup_{u^i\in U^i} f^i(x,u^i) 
\geq \inf_{x\in X} \sum_{t=1}^{T} \theta_t \sum_{i=1}^{m} y^{(i)}_t f^i(x,u_t^i) 
\geq \max_{i\in[m]} \sum_{t=1}^{T} \theta_t f^i(x_t,u_t^i) - (1-\tau) \epsilon 
>0, 
\]
where the first inequality follows from inequality~\eqref{eqn:epsPos}, the second inequality from \eqref{eqn:epsbullet-bound} and the last inequality holds because $\max_{i\in[m]} \sum_{t=1}^{T} \theta_t f^i(x_t,u_t^i) > (1-\tau) \epsilon$. This implies \eqref{eqn:approx-robust-feas} is infeasible.
\Halmos
\endproof

In Section \ref{sec:RegretforRO} we will show that $\epsilon^\circ(\{x_t,u_t,\theta_t\}_{t=1}^T)$ can be interpreted as a weighted regret term \eqref{eqn:weighted-regret-defn}. 
On the other hand, the term $\epsilon^\bullet(\{x_t,u_t,y_t,\theta_t\}_{t=1}^T)$ has no such direct interpretation. In order to upper-bound it by a weighted regret term, we need the following result. 
\begin{corollary}\label{cor:onlineSP-upperbound}
Given sequences $\{x_t,u_t,\theta_t\}_{t=1}^T$ with $x_t\in X$, $u_t\in U$,  for all $t\in[T]$, $\theta\in\Delta_T$, there is an appropriate choice of sequence $\{\bar y_t\}_{t=1}^T$ where $\bar y_t\in\Delta_m$ for all $t\in[T]$, such that $\epsilon^\bullet(\{x_t,u_t,\bar y_t,\theta_t\}_{t=1}^T)$ is upper-bounded by
\begin{equation}\label{eqn:onlineSPgap-upperbound}
\epsilon^\bullet(\{x_t,u_t,\theta_t\}_{t=1}^T) := \sum_{t=1}^{T} \theta_t \max_{i \in [m]} f^i(x_t,u_t) - \inf_{x \in X} \sum_{t=1}^T \theta_t \max_{i \in [m]} f^i(x,u_t).
\end{equation}
Thus, if $\epsilon^\bullet(\{x_t,u_t,\theta_t\}_{t=1}^T)\leq (1-\tau) \epsilon$ and $\max_{i\in[m]} \sum_{t=1}^{T} \theta_t f^i(x_t,u_t^i) > (1-\tau) \epsilon$, then \eqref{eqn:approx-robust-feas} is infeasible.
\end{corollary}
\proof{Proof.}
Given $\{u_t\}_{t=1}^T$, let $x^* \in \argmin_{x \in X} \sum_{t=1}^T \theta_t \max_{i \in [m]} f^i(x,u_t^i)$ and for all $t\in[T]$ define $\bar y_t\in\R^m$ to be the $i$-th unit vector  where $i$ is the smallest index satisfying $i \in \argmax_{i' \in [m]} f^i(x^*,u_t^{i'})$. Then $\inf_{x\in X} \sum_{t=1}^{T} \theta_t \sum_{i=1}^{m} \bar y^{(i)}_t f^i(x,u_t^i) = \inf_{x \in X} \sum_{t=1}^T \theta_t \max_{i \in [m]} f^i(x,u_t^i)$, and the bound follows from 
$\max_{i \in [m]} \sum_{t=1}^{T} \theta_t f^i(x_t,u_t) \leq \sum_{t=1}^{T} \theta_t \max_{i \in [m]} f^i(x_t,u_t)$. We deduce the last result from Theorem~\ref{thm:robust-feas-oracle}.
\Halmos
\endproof

The following corollary demonstrates how we can choose $\tau$ in Theorem~\ref{thm:robust-feas-oracle}.
\begin{corollary}\label{cor:feasCertificates}
Suppose $\{x_t,u_t,y_t,\theta_t\}_{t=1}^T$ with $x_t\in X$, $u_t\in U$, $y_t\in\Delta_m$ for all $t\in[T]$, and $\theta\in\Delta_T$ is such that there exists $\kappa^\circ,\kappa^\bullet\in(0,1)$ satisfying $\epsilon^\circ(\{x_t,u_t,\theta_t\}_{t=1}^T) \leq \epsilon\,\kappa^\circ$ and $\epsilon^\bullet(\{x_t,u_t,y_t,\theta_t\}_{t=1}^T)\leq \epsilon\,\kappa^\bullet$ with $\kappa^\circ+\kappa^\bullet\leq 1$. Let $\tau \in [\kappa^\circ,1-\kappa^\bullet]$. 
Whenever $\max_{i\in[m]} \sum_{t=1}^{T} \theta_t f^i(x_t,u_t^i)\leq (1-\tau) \epsilon$ as well, the solution $\bar{x}_T := \sum_{t=1}^{T} \theta_t x_t$ is $\epsilon$-feasible with respect to \eqref{eqn:approx-robust-feas}. Also, whenever $\max_{i\in[m]} \sum_{t=1}^{T} \theta_t f^i(x_t,u_t^i) > (1-\tau) \epsilon$, then \eqref{eqn:approx-robust-feas} is infeasible.
\end{corollary}
\proof{Proof.}
Note that $\tau \in (0,1)$ follows from its definition, $\kappa^\circ,\kappa^\bullet\geq0$, and $\kappa^\circ+\kappa^\bullet\leq 1$. Furthermore, the interval $[\kappa^\circ,1-\kappa^\bullet]$ is well-defined since $\kappa^\circ \leq 1 - \kappa^\bullet$ always holds. Moreover, $\epsilon^\circ(\{x_t,u_t,\theta_t\}_{t=1}^T) \leq \epsilon\,\kappa^\circ \leq \epsilon\,\tau$ and $\epsilon^\bullet(\{x_t,u_t,y_t,\theta_t\}_{t=1}^T)\leq \epsilon\,\kappa^\bullet \leq \epsilon(1-\tau)$ holds from the definition of $\tau$. The result now follows from Theorem~\ref{thm:robust-feas-oracle}.
\Halmos
\endproof

Theorem~\ref{thm:robust-feas-oracle} and Corollary~\ref{cor:onlineSP-upperbound} points to our general iterative framework for finding robust feasibility/infeasibility certificates of \eqref{eqn:approx-robust-feas}: generate sequences $\{x_t,u_t\}_{t=1}^T$ iteratively to bound $\epsilon^\circ(\{x_t,u_t,\theta_t\}_{t=1}^T)$ and $\epsilon^\bullet(\{x_t,u_t,\theta_t\}_{t=1}^T)$, and then evaluate the term $\max_{i \in [m]} \sum_{t=1}^T \theta_t f^i(x_t,u_t^i)$.
We provide a description of our framework in Algorithm \ref{alg:approx-robust-feas}. We assume that we have access to weights $\{\theta_t\}_{t=1}^T$ and update algorithms $\CA_i$ and $\CA_x$ for choosing $u_t^i \in U^i$ and $x_t \in X$ based on past observations $\{x_s,u_s\}_{s=1}^{t-1}$. We denote the updates by
\[ 
u_t^i = \CA_i(\{x_s,u_s\}_{s=1}^{t-1}) \in U^i~\forall i\in[m], \quad x_t = \CA_x(\{x_s,u_s\}_{s=1}^{t-1}) \in X, \]
and initializations $u_1^i = \CA_i(\{\})~\forall i\in[m]$, $x_1 = \CA_x(\{\})$. Moreover, we assume that these algorithms enjoy the following convergence guarantees: for any sequence $\{x_t\}_{t=1}^T$, let $u_t^i = \CA_i(\{x_s,u_s\}_{s=1}^{t-1})~\forall i\in[m]$, then 
\begin{equation}\label{eqn:noise-update-regret-bound}
\sup_{u^i\in U^i} \sum_{t=1}^T \theta_t f^i(x_t,u^i) - \sum_{t=1}^T \theta_t f^i(x_t,u_t^i) \leq \CR_i(T);
\end{equation}
for any sequence $\{u_s\}_{s=1}^T$, let $x_t = \CA_x(\{x_s,u_s\}_{s=1}^{t-1})$, then
\begin{equation}\label{eqn:solution-update-regret-bound}
\epsilon^\bullet(\{x_t,u_t,\theta_t\}_{t=1}^T) = \sum_{t=1}^T \theta_t \max_{i \in [m]} f^i(x_t,u_t^i) - \inf_{x\in X} \sum_{t=1}^T \theta_t \max_{i \in [m]} f^i(x,u_t^i) \leq \CR_x(T).
\end{equation}
Explicit examples of $\CA_i, \CA_x$ and their bounds $\CR_i, \CR_x$ will be discussed in Section \ref{sec:OCOforRO}. Generally, we desire that the error bounds $\CR_i(T), \CR_x(T) \to 0$ as $T \to \infty$, 
which can be achieved by using online mirror descent as in Theorem~\ref{thm:OCO-non-smooth}. That said, our OFO-based approach in Algorithm \ref{alg:approx-robust-feas} is quite flexible in terms of the selection of OFO algorithms $\CA_i, \CA_x$, and is certainly not restricted to only using online mirror descent.

\begin{algorithm}[ht!]
\caption{OFO-based approximate robust feasibility solver.}\label{alg:approx-robust-feas}
\begin{algorithmic}
{\small
\STATE {\bf input:~} update algorithms $\CA_i$, $i \in [m]$, $\CA_x$, tolerance level $\epsilon>0$, sufficiently large $T = T(\epsilon)$ such that $\max_{i \in [m]} \CR_i(T) + \CR_x(T) \leq \epsilon$, and convex combination weights $\theta_1,\ldots,\theta_T > 0$.
\STATE {\bf output:~} either $\bar{x} \in X$ such that $\sup_{u^i\in U^i} f^i(\bar{x},u^i) \leq \epsilon$ for all $i\in[m]$, or an infeasibility certificate for \eqref{eqn:approx-robust-feas}.
\STATE initialize $u_1^i = \CA_i(\{\})$ for $i\in[m]$ and $x_1 = \CA_x(\{\})$.
\FOR{$t=2,\ldots,T$}
	\FOR{$i=1,\ldots,m$}
		\STATE compute $u_t^i = \CA_i(\{x_s,u_s\}_{s=1}^{t-1}) \in U^i$.
	\ENDFOR
	\STATE compute $x_t = \CA_x(\{x_s,u_s\}_{s=1}^{t-1}) \in X$.
	\STATE obtain upper bounds $\max_{i \in [m]} \CR_i(t) \geq \epsilon^\circ(\{x_s,u_s,\theta_s\}_{s=1}^t)$ and $\CR_x(t) \geq \epsilon^\bullet(\{x_s,u_s,\theta_s\}_{s=1}^t)$.
	\STATE compute $\kappa_t^\circ = \max_{i \in [m]} \CR_i(t)/\epsilon$, $\kappa_t^\bullet = \CR_x(t)/\epsilon$.
	\IF{$\kappa_t^\circ + \kappa_t^\bullet \leq 1$}
		\STATE set $\vartheta_t := \max_{i\in[m]} \sum_{s=1}^{t} \theta_s f^i(x_s,u_s^i)$ and 
		$\tau_t := 1 - \kappa_t^\bullet$.
		\LINEIF{$\vartheta_t > (1-\tau_t)\epsilon$}{\algorithmicreturn~ `infeasible'.}
		\LINEIF{$\vartheta_t \leq (1-\tau_t) \epsilon$}{\algorithmicreturn~ $\bar{x}_t = \frac{1}{t} \sum_{s=1}^t x_s$ as a robust $\epsilon$-feasible solution to \eqref{eqn:approx-robust-feas}.}
	\ENDIF
\ENDFOR
}
\end{algorithmic}
\end{algorithm}

\begin{remark}\label{rem:anticipatory-requirements}
Notice that in Algorithm \ref{alg:approx-robust-feas} we generate $u_t$ \emph{before} generating $x_t$. Thus, nothing stops us from choosing $x_t$ based on the knowledge of $u_t$, or vice versa. Indeed, the pessimization oracle approach of \cite{MutapcicBoyd2009} and the nominal feasibility oracle approach of \cite{BenTalHazan2015} fit within our framework if we rewrite Algorithm \ref{alg:approx-robust-feas} to reflect this, and it is a trivial matter to do so.
However, a conflict may arise if we encounter a situation where generating $x_t$ requires knowledge of $u_t$, and generating $u_t$ also requires knowledge of $x_t$. Thus, when selecting the update algorithms $\CA_i$, $\CA_x$, care must be taken to avoid such situations.
Our suggested OFO approach in Section \ref{sec:RegretforRO} utilizes Theorem~\ref{thm:OCO-non-smooth}. 
Moreover, Theorem~\ref{thm:OCO-non-smooth} generates the current decision in a \emph{non-anticipatory manner} based on only the knowledge of $f_{t-1}$, and not of $f_t$. 
That is, it ensures that we will only use $u_{t-1}$ to generate $x_t$, and similarly we only use $x_{t-1}$ to generate $u_t$, thus no conflicts will arise. 
\epr
\end{remark}

\begin{remark}\label{rem:tauInfeas}
Note that Algorithm~\ref{alg:approx-robust-feas} chooses $\tau_t = 1 - \kappa_t^\bullet$, whereas Corollary~\ref{cor:feasCertificates} allows us to choose from a range $\tau_t \in [\kappa_t^\circ,1-\kappa_t^\bullet]$. This is because it is theoretically possible for \eqref{eqn:approx-robust-feas} to simultaneously be infeasible and robust $\epsilon$-feasible, but in practice we would like to discover infeasibility of \eqref{eqn:approx-robust-feas} rather than an approximately feasible solution. Then the best value for $\tau_t\in[\kappa_t^\circ,1-\kappa_t^\bullet]$ in detecting infeasibility of \eqref{eqn:approx-robust-feas} is given by $\tau_t = 1 - \kappa_t^\bullet$.
\epr
\end{remark}

In the next section, we describe some approaches to implement Algorithm \ref{alg:approx-robust-feas} in practice by providing explicit examples of $\CA_i$ and $\CA_x$.

\section{Customizations of the General Framework}\label{sec:OCOforRO}

In this section, we examine how to generate the sequences $\{x_t,u_t\}_{t=1}^T$ in practice. In Section~\ref{sec:RegretforRO}, we first interpret the terms $\epsilon^\circ(\{x_t,u_t,\theta_t\}_{t=1}^T)$ and $\epsilon^\bullet(\{x_t,u_t,\theta_t\}_{t=1}^T)$ from Section~\ref{sec:SPFormulation} as weighted regret terms, which gives rise to our OFO-based approach. In Section \ref{sec:pessimization-oracle}, we modify the pessimization oracle-based approach of \cite{MutapcicBoyd2009} to solving \eqref{eqn:approx-robust-feas} within our framework. In Section \ref{sec:oracle-based-approach}, we examine the nominal feasibility oracle-based approach of \cite{BenTalHazan2015} within the context of our general framework. Finally, in Section~\ref{sec:RateDiscussion}, we summarize and compare the convergence rates achievable via various customizations of our framework using these different approaches.

\subsection{The OFO-based Approach}\label{sec:RegretforRO}
Let us first consider $\epsilon^\circ(\{x_t,u_t,\theta_t\}_{t=1}^T)$. For any $i\in[m]$, given $x_t$, we define the function $f_t^i: U^i\to \R$ as $f_t^i(u^i) = -f^i(x_t,u^i)$. Then the function $f_t^i(u^i)$ is convex in $u^i$ under Assumption~\ref{ass:f-concave-u}, and the subterm of $\epsilon^\circ(\{x_t,u_t,\theta_t\}_{t=1}^T)$ given by 
\begin{equation}\label{eqn:RO-weighted-regret}
\sup_{u^i\in U^i} \sum_{t=1}^{T} \theta_t f^i(x_t,u^i) - \sum_{t=1}^{T} \theta_t f^i(x_t,u_t^i)
\end{equation}
is the weighted regret \eqref{eqn:weighted-regret-defn} corresponding to the sequence of functions $\{ f_t^i \}_{t=1}^T$. When the uncertainty sets $U^i$, $i \in [m]$ admit proximal setups as in Assumption~\ref{ass:proximal-setup-OFO}, Theorem~\ref{thm:OCO-non-smooth} from Section~\ref{sec:OCO} gives an efficient OFO algorithm for choosing $\{u_t^i\}_{i=1}^m$ to bound the regret subterms \eqref{eqn:RO-weighted-regret} with $O(1/\sqrt{T})$. 
Therefore, by using the online mirror descent algorithm of Theorem \ref{thm:OCO-non-smooth} as $\CA_i$ in the computation of our $u^i_t$, we guarantee that
\[ \epsilon^\circ(\{x_t,u_t,\theta_t\}_{t=1}^T) = \max_{i \in [m]} \left\{ \sup_{u^i\in U^i} \sum_{t=1}^{T} \theta_t f^i(x_t,u^i) - \sum_{t=1}^{T} \theta_t f^i(x_t,u_t^i) \right\} \leq \max_{i \in [m]} \CR_i(T), \]
where $\CR_i(T) = O(1/\sqrt{T})$ with uniform weights $\theta_t = 1/T$.

On the other hand, given $u_t^i\in U^i$ for $i\in[m]$, let us define $\varphi_t(x) := \max_{i \in [m]} f^i(x,u_t^i)$. Then $\varphi_t(x)$ is convex in $x$ over $X$ since the functions $f^i$ are convex in $x$ by Assumption~\ref{ass:f-concave-u}. We can then rewrite $\epsilon^\bullet(\{x_t,u_t,\theta_t\}_{t=1}^T)$ as
\begin{align}\label{eqn:RO-onlineSPgap-upperbd}
\epsilon^\bullet(\{x_t,u_t,\theta_t\}_{t=1}^T) &= \sum_{t=1}^{T} \theta_t \max_{i \in [m]} f^i(x_t,u_t^i) - \inf_{x\in X} \sum_{t=1}^{T} \theta_t \max_{i \in [m]} f^i(x,u_t^i)\notag\\
 &= \sum_{t=1}^{T} \theta_t \varphi_t(x_t) - \inf_{x\in X} \sum_{t=1}^{T} \theta_t \varphi_t(x).
\end{align}
Then $\epsilon^\bullet(\{x_t,u_t,\theta_t\}_{t=1}^T)$ is also a weighted regret term \eqref{eqn:weighted-regret-defn} corresponding to the sequence of functions $\{\varphi_t\}_{t=1}^T$. When the domain $X$ admits a proximal setup as in Assumption~\ref{ass:proximal-setup-OFO}, Theorem~\ref{thm:OCO-non-smooth} again gives an efficient OFO algorithm for choosing $x_t$ to bound \eqref{eqn:RO-onlineSPgap-upperbd}. 
Once again, we may choose $\CA_x$ to be the online mirror descent, and get $\CR_x(T) = O(1/\sqrt{T})$ with uniform weights $\theta_t = 1/T$.

Algorithm \ref{alg:approx-robust-feas} can now be employed, provided we choose $T = \Omega(1/\epsilon^2)$, to solve the robust feasibility problem \eqref{eqn:approx-robust-feas}. Since the online mirror descent algorithm of Theorem~\ref{thm:robust-feas-oracle} only uses first-order information, we can solve the robust feasibility problem \eqref{eqn:approx-robust-feas} while avoiding reliance on a pessimization oracle for $u$ as in \cite{MutapcicBoyd2009} or a nominal feasibility oracle for $x$ as in \cite{BenTalHazan2015}.

\subsection{The Pessimization Oracle-Based Approach}\label{sec:pessimization-oracle}
Mutapcic and Boyd \cite{MutapcicBoyd2009} generate solutions $x_t \in X$ at each iteration $t$ by solving an extended nominal problem
\begin{equation}\label{eqn:extended-nominal-solver}
\min_{x \in X} \left\{ f^0(x) :~ f^i(x,u^i) \leq 0, \ \forall u^i \in \hat{U}_{t-1}^i, \  i \in [m] \right\},
\end{equation}
where $\hat{U}_{t-1}^i \subset U^i$ are finite approximate uncertainty sets based on past noise realizations $\{u_{t'}^i\}_{i=1}^m$ for $t'\in[t-1]$. New noises $u_t$ are then generated by calling the \emph{pessimization oracles} on the current solution $x_t$. More precisely, given $x_t \in X$, the pessimization oracles solve $\sup_{u^i \in U^i} f^i(x_t,u^i)$ and return
\begin{equation}\label{eqn:pessimization-oracle}
u_t^i \in U^i \quad \text{s.t.} \quad f^i(x_t,u_t^i) \geq \sup_{u^i \in U^i} f^i(x_t,u^i) - \tau \epsilon.
\end{equation}
In terms of our framework of Algorithm~\ref{alg:approx-robust-feas}, the update policy of generating new noises $u_t$ in this approach of \cite{MutapcicBoyd2009} corresponds to selecting the algorithms $\CA_x$ to be a extended nominal solver for \eqref{eqn:extended-nominal-solver} and the algorithms $\CA_i$ to be pessimization oracles that solve \eqref{eqn:pessimization-oracle}. Note that computing $u_t^i$ requires knowledge of $x_t$ (see Remark \ref{rem:anticipatory-requirements}), and consequently the bound for the regret term \eqref{eqn:noise-update-regret-bound} is $\CR_i(T) \leq \tau \epsilon$ for any $T$. We show this in the proof of Theorem \ref{thm:pessimization-oracle-approach}. If for all $i \in [m]$ we have $f^i(x_t,u_t^i) \leq (1-\tau) \epsilon$, then we terminate and declare $x_t$ is a robust $\epsilon$-feasible and optimal solution; otherwise, we append $\hat{U}_t^i = \hat{U}_{t-1}^i \cup \{u_t^i\}$ and re-solve \eqref{eqn:extended-nominal-solver} with the new approximate sets $\hat{U}_t^i$. It is shown in \cite[Section 5.2]{MutapcicBoyd2009} that the number of iterations $T$ needed before termination with a robust $\epsilon$-feasible solution $x_T$ is upper bounded by $(1 + O(1/\epsilon))^n$ where $n$ is the dimension of $x$.

Suppose now that we are interested in robust feasibility \eqref{eqn:approx-robust-feas}. \cite[Section 5.3]{MutapcicBoyd2009} discusses a number of variations for generating $x_t$ by modifying \eqref{eqn:extended-nominal-solver}. In contrast, we propose the following modification: instead of solving \eqref{eqn:extended-nominal-solver}, generate $\{x_t\}_{t=1}^T$ via a non-anticipatory algorithm $\CA_x$ (see Remark \ref{rem:anticipatory-requirements}) to bound $\epsilon^\bullet(\{x_t,u_t,\theta_t\}_{t=1}^T) \leq \CR_x(T) = (1-\tau)\epsilon$. We call  our modification \emph{FO-based pessimization}. Then the FO-based pessimization approach fits within our framework as a special case. 
\begin{theorem}\label{thm:pessimization-oracle-approach}
Let $\tau \in (0,1)$. Suppose $\{x_t\}_{t=1}^T$ are generated iteratively to guarantee that $\epsilon^\bullet(\{x_t,u_t,\theta_t\}_{t=1}^T) \leq (1-\tau)\epsilon$ for \emph{any} sequence $\{u_t\}_{t=1}^T$. Suppose $u_t^i$ are generated by pessimization oracles \eqref{eqn:pessimization-oracle} for $i \in [m]$. If there exists $t \in [T]$ such that for all $i \in [m]$ we have $f^i(x_t,u_t^i) \leq (1-\tau) \epsilon$, then $x_t$ is a robust $\epsilon$-feasible solution to \eqref{eqn:approx-robust-feas}. If $\max_{i \in [m]} \sum_{t=1}^{T} \theta_t f^i(x_t,u_t^i) \leq (1-\tau) \epsilon$, then $\bar{x}_T = \sum_{t=1}^{T} \theta_t x_t$ is a robust $\epsilon$-feasible solution to \eqref{eqn:approx-robust-feas}. If $\max_{i \in [m]} \sum_{t=1}^{T} \theta_t f^i(x_t,u_t^i) > (1-\tau) \epsilon$, then we certify that \eqref{eqn:approx-robust-feas} is robust infeasible.
\end{theorem}
\proof{Proof.}
It is clear that if there exists $t \in [T]$ such that for all $i \in [m]$ we have $f^i(x_t,u_t^i) \leq (1-\tau) \epsilon$, then $x_t$ is a robust $\epsilon$-feasible solution to \eqref{eqn:approx-robust-feas}. Furthermore, the fact that $\max_{i \in [m]} \sum_{t=1}^{T} \theta_t f^i(x_t,u_t^i) > (1-\tau) \epsilon$ implies robust infeasibility of \eqref{eqn:approx-robust-feas} follows from our assumption that $\epsilon^\bullet(\{x_t,u_t,\theta_t\}_{t=1}^T) \leq (1-\tau) \epsilon$ and Corollary \ref{cor:onlineSP-upperbound}. To show that $\max_{i \in [m]} \sum_{t=1}^{T} f^i(x_t,u_t^i) \leq (1-\tau) \epsilon$ implies that $\bar{x}_T$ is robust $\epsilon$-feasible, we only need to show that $\epsilon^\circ(\{x_t,u_t,\theta_t\}_{t=1}^T) \leq \tau \epsilon$. Observe that by our definition of $u_t^i$ in \eqref{eqn:pessimization-oracle}, we have $f^i(x_t,u_t^i) \geq \sup_{u^i \in U^i} f^i(x_t,u^i) - \tau \epsilon$, hence the regret terms in $\epsilon^\circ(\{x_t,u_t,\theta_t\}_{t=1}^T)$ satisfy
\begin{align*}
\sup_{u^i \in U^i} \sum_{t=1}^{T} \theta_t f^i(x_t,u^i) - \sum_{t=1}^{T} \theta_t f^i(x_t,u_t^i) &\leq \sup_{u^i \in U^i} \sum_{t=1}^{T} \theta_t f^i(x_t,u^i) - \sum_{t=1}^{T} \theta_t \left(\sup_{u^i \in U^i} f^i(x_t,u^i) - \tau \epsilon\right)\\
&= \sup_{u^i \in U^i} \sum_{t=1}^{T} \theta_t f^i(x_t,u^i) - \sum_{t=1}^{T} \theta_t \sup_{u^i \in U^i} f^i(x_t,u^i) + \tau \epsilon\\
&\leq \tau \epsilon.
\end{align*}
Then 
\[ \epsilon^\circ(\{x_t,u_t,\theta_t\}_{t=1}^T) = \max_{i \in [m]} \left\{ \sup_{u^i \in U^i} \sum_{t=1}^{T} \theta_t f^i(x_t,u^i) - \sum_{t=1}^{T} \theta_t f^i(x_t,u_t^i) \right\} \leq \tau \epsilon, \]
and the result follows from Corollary \ref{cor:onlineSP-upperbound}.
\Halmos
\endproof

Theorem~\ref{thm:pessimization-oracle-approach} can only be used to certify robust feasibility/infeasibility. Hence, to find a robust $\epsilon$-optimal solution in FO-based pessimization approach, we must perform a binary search and solve at most $O(\log(1/\epsilon))$ instances of robust feasibility problems. Despite this, in Section~\ref{sec:RateDiscussion}, we discuss how FO-based pessimization approach which uses  OFO algorithms to generate $\{x_t\}_{t=1}^T$ to bound $\epsilon^\bullet(\{x_t,u_t,\theta_t\}_{t=1}^T)$ results in much better complexity guarantees than using an extended nominal feasibility solver  \eqref{eqn:extended-nominal-solver} as proposed by \cite{MutapcicBoyd2009}, even when taking into account the additional $O(\log(1/\epsilon))$ factor. 

\begin{remark}\label{rem:anticipatory-pessimization-oracle}
In the pessimization oracle-based approach, the noises $u_t$ need to be generated with knowledge of $x_t$, because it is not possible to guarantee $f^i(x_t,u_t^i) \geq \sup_{u^i \in U^i} f^i(x_t,u^i) - \tau \epsilon$ if the vectors $u_t^i$ were chosen with only the knowledge of $x_1,\ldots,x_{t-1}$.
\epr
\end{remark}

\subsection{The Nominal Feasibility Oracle-Based Approach}\label{sec:oracle-based-approach}

The nominal feasibility oracle-based approach of Ben-Tal et al.\@ \cite{BenTalHazan2015} suggest using OFO algorithms to choose a sequence $\{u_t\}_{t=1}^T$ that guarantees $\epsilon^\circ(\{x_t,u_t,\theta_t\}_{t=1}^T)$ is small, in a non-anticipatory fashion, for \emph{any} sequence $\{x_t\}_{t=1}^T$. In this aspect, it essentially matches with our OFO-based approach outlined in Section~\ref{sec:RegretforRO} i.e., the choice of $\CA_i$ is essentially the same. The key differentiating point between our OFO-based approach and that of \cite{BenTalHazan2015} lies in which algorithm is chosen for $\CA_x$. At step $t$, \cite{BenTalHazan2015} utilizes a \emph{nominal feasibility oracle}. That is, given parameters $u_t$, they call a powerful, and potentially expensive, nominal feasibility oracle that solves the following feasibility problem to $\epsilon$-accuracy
\begin{equation}\label{eqn:nominal-feas-oracle}
\begin{cases}
\text{\emph{Either}: find}\ \ x \in X \quad\text{s.t.}\quad f^i(x,u_t^i) \leq (1-\tau) \epsilon \quad \forall i\in[m];\\
\text{\emph{or}: declare infeasibility, } \forall x \in X,\ \exists i \in [m] \quad \text{s.t.} \quad f^i(x,u_t^i) > 0.
\end{cases} 
\end{equation}
We denote $x_t \in X$ to be the point returned by this oracle at step $t$, if it exists. For this approach, the outputs of a nominal feasibility oracle can be used to deduce a result similar to Corollary \ref{cor:onlineSP-upperbound}, except that we no longer need to evaluate $\epsilon^\bullet(\{x_t,u_t,\theta_t\}_{t=1}^T)$, we just need to bound $\epsilon^\circ(\{x_t,u_t,\theta_t\}_{t=1}^T)$.
\begin{theorem}\label{thm:robust-feas-oracle-solver}
Given weights $\theta\in\Delta_T$, suppose that the sequence $\{u_t\}_{t=1}^T$ is generated in a non-anticipatory manner to guarantee $\epsilon^\circ(\{x_t,u_t,\theta_t\}_{t=1}^T) \leq \tau \epsilon$ for any sequence $\{x_t\}_{t=1}^T$. Also, suppose that at each step $t \in [T]$, $x_t$ is generated by the nominal feasibility oracle which solves \eqref{eqn:nominal-feas-oracle}. If there exists $t \in [T]$ such that \eqref{eqn:nominal-feas-oracle} declares infeasibility, then \eqref{eqn:approx-robust-feas} is infeasible. Otherwise, if $x_t$ satisfies $f^i(x_t,u_t^i) \leq (1-\tau)\epsilon$ for all $t \in [T]$ and $i \in [m]$, we have a robust $\epsilon$-feasibility certificate for \eqref{eqn:approx-robust-feas}.
\end{theorem}
\proof{Proof.}
If \eqref{eqn:nominal-feas-oracle} declares infeasibility, then it is obvious that the robust feasibility problem is infeasible. We focus on the latter case. By the premise of the theorem, we have $\epsilon^\circ(\{x_t,u_t,\theta_t\}_{t=1}^T) \leq \tau\epsilon$. Let us evaluate $\max_{i\in[m]} \sum_{t=1}^{T} \theta_t f^i(x_t,u_t^i)$. Because $\theta \in \Delta_T$ and from the definition of the nominal feasibility oracle we have $f^i(x_t,u_t^i) \leq (1-\tau) \epsilon$ for all $t \in [T]$ and $i \in [m]$, we conclude $\max_{i\in[m]} \sum_{t=1}^{T} \theta_t f^i(x_t,u_t^i) \leq (1-\tau) \epsilon$. The conclusion now follows from Theorem~\ref{thm:robust-feas-oracle}.
\Halmos
\endproof 
Thus, the approach of \cite{BenTalHazan2015}, which works with nominal feasibility oracles, fits within our framework of Algorithm~\ref{alg:approx-robust-feas} right away. We next make three important remarks.

\begin{remark}\label{rem:anticipatory-nom-feas-oracle}
Similar to Remark~\ref{rem:anticipatory-pessimization-oracle}, a critical property required in the approach of \cite{BenTalHazan2015} of the vectors $x_t$ is that $f^i(x_t,u_t^i) \leq (1-\tau) \epsilon$. This is possible only if $x_t$ were chosen with the knowledge of $\{u_1^i\}_{i=1}^m,\ldots,\{u_{t}^i\}_{i=1}^m$. 
\epr
\end{remark}

\begin{remark}\label{rem:oracle-based-optimal-sol}
Theorem~\ref{thm:robust-feas-oracle-solver} states that the nominal feasibility oracle-based approach can solve robust feasibility problems \eqref{eqn:approx-robust-feas}. This then recovers \cite[Theorems 1,2]{BenTalHazan2015}. In addition, we next make a nice and practical  observation that was overlooked in \cite{BenTalHazan2015}. We show that slightly adjusting this oracle will let us directly solve the robust \emph{optimization} problem \eqref{eqn:robust-opt-intro}, i.e., optimize a convex objective function $f^0(x)$ instead of relying on a binary search over the optimal objective value. Recall that $\Opt$ is the optimal value of the RO problem (see \eqref{eqn:robust-opt-intro}). Naively, to solve for $\Opt$, we would embed $f^0$ into the constraint set, and then perform a binary search over the robust feasible set by repeatedly applying the oracle-based approach and Theorem~\ref{thm:robust-feas-oracle-solver} to check for robust feasibility. Suppose  that now, instead of using a nominal feasibility oracle to solve \eqref{eqn:nominal-feas-oracle}, we work with a \emph{nominal optimization oracle}. That is, given fixed parameters $\{u_t^i\}_{i=1}^m$, we have access to an oracle that solves
\[ \Opt_t = \inf_x \left\{ f^0(x) : f^i(x,u_t^i) \leq 0, \ i \in [m],\ \ x \in X \right\}. \]
When solving for $\Opt_t$, most convex optimization solvers will either declare that the constraints are infeasible, or return a point $x_t \in X$ such that $f^i(x_t,u_t^i) \leq (1-\tau)\epsilon$ and $f^0(x_t) \leq \Opt_t + \epsilon$. It is clear that $f^0(x_t) \leq \Opt_t + \epsilon \leq \Opt + \epsilon$. Given such a sequence of points $\{x_t\}_{t=1}^T$, from Theorem~\ref{thm:robust-feas-oracle-solver} we deduced that $\bar{x}_T = \sum_{t=1}^{T} \theta_t x_t$ is a robust $\epsilon$-feasible solution. Moreover,  convexity of $f^0$ implies 
\[ f^0(\bar{x}_T) \leq \sum_{t=1}^{T} \theta_t f^0(x_t) \leq \sum_{t=1}^{T} \theta_t (\Opt + \epsilon) = \Opt + \epsilon. \]
Hence, not only do we claim that $\bar{x}_T$ is robust $\epsilon$-feasible, but that it is also $\epsilon$-optimal. Thus, when our oracle can return $\epsilon$-optimal solutions, which most solvers can, we eliminate the need to perform a binary search.
\epr
\end{remark}

Below we elaborate on the differences between Theorem~\ref{thm:robust-feas-oracle-solver} and Corollary \ref{cor:onlineSP-upperbound}. 
\begin{remark}\label{rem:nominal-feas-oracle-inefficient-SV}
In contrast to Corollary \ref{cor:onlineSP-upperbound}, Theorem~\ref{thm:robust-feas-oracle-solver} does not need to control the term $\epsilon^\bullet(\{x_t,u_t,\theta_t\}_{t=1}^T)$. There are two reasons for this: $(i)$ due to \eqref{eqn:nominal-feas-oracle}, each point $x_t$ satisfies $f^i(x_t,u_t^i) \leq (1-\tau)\epsilon$, hence $\max_{i\in[m]} \sum_{t=1}^{T} \theta_t f^i(x_t,u_t^i) \leq (1-\tau)\epsilon$ always holds. Therefore, the infeasibility part of Corollary \ref{cor:onlineSP-upperbound} never becomes relevant, and $(ii)$ due to the oracle solving \eqref{eqn:nominal-feas-oracle}, infeasibility may be declared at any step $t \in [T]$ in Theorem~\ref{thm:robust-feas-oracle-solver}. This offers the possibility of stopping  early rather than having to wait until all $T$ steps are completed. 
Thus, the nominal feasibility oracle-based approach trades off using more effort at each iteration $t$ to solve \eqref{eqn:nominal-feas-oracle} for the ability to terminate early. In contrast, our OFO-based approach opts to keep the per-iteration cost cheap while giving up the ability to terminate early. 
More formally, let us examine a particular way of  solving \eqref{eqn:nominal-feas-oracle} within a nominal feasibility oracle. Note that \eqref{eqn:nominal-feas-oracle} is equivalent to checking $F_t\leq (1-\tau) \epsilon$ or $F_t>0$, where
\begin{equation}\label{eqn:SV-feas-points}
F_t := \inf_{x \in X} \left\{ \max_{i \in [m]} f^i(x,u_t^i) \right\}.
\end{equation}
Since each $f^i(x,u_t^i)$ is convex in $x$ for fixed $u_t^i$, $\max_{i \in [m]} f^i(x,u_t^i)$ is convex in $x$ also, hence standard convex optimization methods may be employed to find $x_t \in X$ such that
\[ F_t \leq \max_{i \in [m]} f^i(x_t,u_t^i) \leq F_t + (1-\tau) \epsilon. \]
Then, by checking whether $\max_{i \in [m]} f^i(x_t,u_t^i) \leq (1-\tau) \epsilon$ or $\max_{i \in [m]} f^i(x_t,u_t^i) > (1-\tau) \epsilon$, we can determine whether $F_t \leq (1-\tau) \epsilon$ or $F_t > 0$ respectively. In particular, if we find that $F_t \leq (1-\tau) \epsilon$, our point $x_t$ is feasible for \eqref{eqn:nominal-feas-oracle}.

Also, when all the vectors $x_t$ satisfy \eqref{eqn:SV-feas-points}, we have the bound
\begin{align*}
\epsilon^\bullet(\{x_t,u_t,\theta_t\}_{t=1}^T) &= \sum_{t=1}^{T} \theta_t \max_{i \in [m]} f^i(x_t,u_t^i) - \inf_{x \in X} \sum_{t=1}^{T} \theta_t \max_{i \in [m]} f^i(x,u_t^i)\\
&\leq \sum_{t=1}^{T} \theta_t \left[ \max_{i \in [m]} f^i(x_t,u_t^i) - \inf_{x \in X} \max_{i \in [m]} f^i(x,u_t^i) \right] \leq (1-\tau) \epsilon.
\end{align*}
Consequently, we deduce that the nominal feasibility oracle, implemented as a convex optimization problem, also naturally bounds $\epsilon^\bullet(\{x_t,u_t,\theta_t\}_{t=1}^T)$ although this bound is not utilized in Theorem~\ref{thm:robust-feas-oracle-solver}. In terms of our framework of Algorithm~\ref{alg:approx-robust-feas}, the update policy of generating new solutions $x_t$ in this approach corresponds to selecting the algorithm $\CA_x$ to be a convex optimization solver that solves $\Opt_t$. Then whenever the solver returns a feasible solution, the regret bound \eqref{eqn:solution-update-regret-bound} satisfies  $\epsilon^\bullet(\{x_t,u_t,\theta_t\}_{t=1}^T) \leq \CR_x(T) = (1-\tau) \epsilon$ for any $T$. Note that the term $\epsilon^\bullet(\{x_t,u_t,\theta_t\}_{t=1}^T)$ inherently includes the objective functions $\max_{i \in [m]} f^i(x,u_t^i)$ of each problem $F_t$. At each iteration $t$, instead of evaluating $F_t$ to $(1-\tau)\epsilon$ accuracy, our OFO-based approach performs only a simple update based on the first-order information, and it yields a bound on $\epsilon^\bullet(\{x_t,u_t,\theta_t\}_{t=1}^T)$ from the overall collection of these simple updates.
\epr
\end{remark}

\subsection{Convergence Rates and Discussion}\label{sec:RateDiscussion}

We summarize the convergence rates achievable in our general RO framework for various cases. 
We first examine the number of iterations required for each approach discussed, then proceed to analyze the per-iteration cost of each approach. A summary of our discussion is given in Table \ref{tab:rate-comparison}. 
We use the notation $r_u(\epsilon)$ to denote the number of iterations
$T$ required for algorithms $\CA_i$ to guarantee $\epsilon^\circ(\{x_t,u_t,\theta_t\}_{t=1}^T) \leq \max_{i \in [m]} \CR_i(T) \leq \epsilon/2$. Similarly, we let $r_x(\epsilon)$ be the number of iterations $T$ required for algorithm $\CA_x$ to guarantee that $\epsilon^\bullet(\{x_t,u_t,\theta_t\}_{t=1}^T) \leq \CR_x(T) \leq \epsilon/2$. Then the resulting worst-case number of iterations needed in Algorithm~\ref{alg:approx-robust-feas} to obtain robust $\epsilon$-feasibility/infeasibility certificates is $\max\{r_u(\epsilon),r_x(\epsilon)\}$.

As outlined in Section \ref{sec:OCOforRO}, employing standard OFO-based algorithms, i.e., Theorem \ref{thm:OCO-non-smooth}, on the terms $\epsilon^\circ(\{x_t,u_t,\theta_t\}_{t=1}^T)$ and $\epsilon^\bullet(\{x_t,u_t,\theta_t\}_{t=1}^T)$ requires $r_u(\epsilon) = O(1/\epsilon^2)$ and $r_x(\epsilon) = O(1/\epsilon^2)$ iterations to ensure they are no larger than $\epsilon/2$. Thus, our OFO-based approach from Section \ref{sec:RegretforRO} requires $O(1/\epsilon^2)$ iterations to solve \eqref{eqn:robust-feas}. Since our OFO-based approach returns only robust $\epsilon$-feasible solutions, we need to perform a binary search and repeatedly invoke our method $O(\log(1/\epsilon))$ times to obtain $\epsilon$-optimal solutions, so the total number of iterations is $O(\log(1/\epsilon)/\epsilon^2)$.

Our FO-based pessimization approach, i.e., our modification of the pessimization oracle-based approach of \cite{MutapcicBoyd2009} outlined in Section \ref{sec:pessimization-oracle}, requires $r_x(\epsilon)$ iterations to solve \eqref{eqn:robust-feas} because by Theorem~\ref{thm:pessimization-oracle-approach} we only need to guarantee $\epsilon^\bullet(\{x_t,u_t,\theta_t\}_{t=1}^T) \leq \CR_x(T) \leq \epsilon/2$. Taking into account the binary search factor $O(\log(1/\epsilon))$ to find a robust $\epsilon$-optimal solution, the total number of iterations required is $O(\log(1/\epsilon)/\epsilon^2)$, which is much better than the exponential $(1 + O(1/\epsilon))^n$ bound of \cite[Section 5.2]{MutapcicBoyd2009} that uses a full nominal solution oracle \eqref{eqn:extended-nominal-solver}. Similarly, the nominal feasibility/optimization oracle-based approach of \cite{BenTalHazan2015} outlined in Section \ref{sec:oracle-based-approach} requires $r_u(\epsilon) = O(1/\epsilon^2)$ iterations (or $r_u(\epsilon)\log(1/\epsilon) = O(\log(1/\epsilon)/\epsilon^2)$ iterations if only a feasibility oracle is used) to obtain robust $\epsilon$-optimal solutions because by Theorem \ref{thm:robust-feas-oracle-solver} we only need to bound $\epsilon^\circ(\{x_t,u_t,\theta_t\}_{t=1}^T) \leq \max_{i \in [m]} \CR_i(T) \leq \epsilon/2$.

\begin{remark}\label{rem:FOMflexibility}
The flexibility of our general framework in terms of the selection of algorithms $\CA
_i,\cA_x$ extends beyond just using Theorem \ref{thm:OCO-non-smooth}. Depending on the structure of functions $f^i$ and uncertainty domains $U^i$, the algorithms $\cA_i$ and $\cA_x$ may be replaced by more appropriate OCO algorithms. For example, when $f^i$ are strongly convex, certain OCO algorithms achieve faster convergence rates. Moreover, unless explicitly required by the algorithms $\cA_i$, we do not need to assume convexity of the sets $U^i$. As a result, the follow-the-leader or follow-the-perturbed-leader type algorithms from \cite{KalaiVempala2005} can be utilized as $\cA_i$ in our framework even when $U^i$ are nonconvex but certain assumptions ensuring applicability of these algorithms are satisfied. Such assumptions are satisfied for example when $f^i(x,u^i)$ are linear in $u^i$ and the nonconvex sets $U^i$ admit a certain linear optimization oracle. This is for example the case in a certain lifted representation of the robust convex quadratic constraint discussed in \cite[Section 4.2]{BenTalHazan2015}.
Similarly, when the functions $f^i(x,u^i)$ are exp-concave in $u^i$, applying the online Newton step algorithm of \cite{HazanAgarwalKale2007} for $\cA_i$ results in a weighted regret bound of at most $O\left( \log(T)/ T \right)$ in $T$ iterations. Such $f^i$ that are exp-concave in $u^i$ satisfying Assumption~\ref{ass:f-concave-u} arise in optimization under uncertainty problems where variance is used as a risk measure, e.g., mean-variance portfolio optimization problems, see for example \cite[Example 25]{BenTalDenHertog2015}. Essentially, the same flexibility for acceleration and/or working with nonconvex sets $U^i$ is present in \cite{BenTalHazan2015} as well. 

In the presence of favorable problem structure, based on Table \ref{tab:rate-comparison}, if an accelerated algorithm to exploit problem structure is employed in the place of $\CA_i$, the overall number of iterations of the nominal feasibility approach is immediately reduced accordingly. Analogous result holds for $\CA_x$ and the FO-based pessimization approach. However, in the case of our OFO-based approach, we need to have favorable structure in \emph{both} $x$ and $u$ and utilize the corresponding accelerated algorithms $\CA_x, \CA_i$ to attain the acceleration of the overall approach.
\epr
\end{remark}

We now discuss the per-iteration cost for each approach. In order to discuss the total \emph{arithmetic complexity} of each approach, we  let $k$ be the maximum dimension of the uncertain parameters $u^i$ for $i\in[m]$ and recall  that $n$ denotes the dimension of the decision variables $x$. In the case where our domains $X,\{U^i\}_{i=1}^m$ have favorable geometry, such as Euclidean ball or simplex, the vectors $x_t, \{u_t^i\}_{i=1}^m$ are updated via simple closed-form prox operations, which cost $O(n)$ and $O(km)$ per iteration respectively. The cost of computing the subgradients $\grad_x f^i(x,u^i), \grad_u f^i(x,u^i)$ is at least $O(km + mn)$ each iteration. This cost is incurred in each iteration of all of the approaches we discuss. From this, we deduce that the per-iteration cost of our OFO-based approach is at most $O(km + mn)$.

The per-iteration cost of the pessimization oracle based approaches involve calling $m$ pessimization oracles \eqref{eqn:pessimization-oracle} and the costs related to updating $x_t$. 
We denote by $\Pess(\epsilon,k)$ the complexity of a pessimization oracle with tolerance $\epsilon$ and $k$ variables. A summary of different possible implementations is given in Table \ref{tab:pess-complexity}.  
If $\sup_{u^i \in U^i} f^i(x,u^i)$ has a simple closed form solution, then the resulting arithmetic cost for $\Pess(\epsilon,k)$ is $O(k)$ for each pessimization oracle. If we can use polynomial-time IPMs, this cost becomes $O(k^3 \log(1/\epsilon))$ (see \cite[Section 6.6]{BenTalNemirovski2001book}), and using FOMs has cost $O(k \log(1/\epsilon))$ in the best case when the functions $f^i$ are smooth \emph{and} strongly convex in  $u^i$.
 In the case of our FO-based pessimization approach, the update involving $x_t$ will be given by simple closed form formulas for prox operations when $X$ has favorable geometry, resulting in a cost of $O(mn)$. The full pessimization approach of \cite{MutapcicBoyd2009} incurs the cost of solving an extended nominal feasibility problem for the update of $x_t$. 

The per-iteration cost of the nominal feasibility/optimization oracle-based approach of \cite{BenTalHazan2015}, as well as that of of \cite{MutapcicBoyd2009}, depends on the type of solver used to solve the nominal optimization/feasibility problem \eqref{eqn:nominal-feas-oracle}. 
We denote by $\Nom(\epsilon,m,n)$ the complexity of a nominal oracle with tolerance $\epsilon$, $m$ constraints and $n$ variables. Note that nominal solvers can be either optimization or feasibility solvers. If it is the latter, an extra $\log(1/\epsilon)$ factor is incurred to perform binary search. A summary of different possible implementations for $\Nom(\epsilon,m,n)$ is given in Table \ref{tab:nom-complexity}. When applicable for $\Nom(\epsilon,m,n)$ implementation, polynomial-time IPMs are guaranteed to terminate in $O(\sqrt{m}\log(1/\epsilon))$ iterations with a solution to \eqref{eqn:nominal-feas-oracle} and thus offer the best rates in terms of their dependence on $\epsilon$. They also have the advantage that they can act as a nominal optimization oracle, and hence by Remark~\ref{rem:oracle-based-optimal-sol} there will be no need to perform an additional binary search to find an $\epsilon$-optimal solution. On the other hand, they demand significantly more memory, and their per-iteration cost is quite high in terms of the dimension, usually around the order of $O(n^3 + mn)$, see \cite[Chapter 6.6]{BenTalNemirovski2001book}. 
In order to keep both the memory requirements and the per-iteration cost associated with implementing the nominal feasibility oracle $\Nom(\epsilon,m,n)$ low, one may opt for a FOM called the CoMirror algorithm that can work with functional constraints, see \cite{BeckBenTal2010} and \cite[Section 1.3]{JuditNem2012Pt1}. CoMirror algorithm is guaranteed to find a solution to the nominal $\epsilon$-feasibility problem within $O(1/\epsilon^2)$ iterations, with a much cheaper per-iteration cost of $O(mn)$. Because CoMirror method can optimize as well, it does not need binary search. However, to the best of our knowledge, its possibility to exploit further structural properties of the functions $f^i$, such as smoothness in $x$, to improve the dependence on $\epsilon$ are not known. In order to exploit such properties in the implementation of $\Nom(\epsilon,m,n)$, it is possible to cast \eqref{eqn:nominal-feas-oracle} as a convex-concave SP problem, and then apply efficient FOMs such as Nesterov's algorithm~\cite{Nesterov2005} or Nemirovski's Mirror Prox algorithm~\cite{Nemirovski2005} to achieve a convergence rate of $O(\log(m)/\epsilon)$ and per-iteration cost of $O(mn)$. This convex-concave SP approach can only be used as a nominal feasibility oracle, so we must repeat the process $\log(1/\epsilon)$ times to obtain an $\epsilon$-optimal solution.

\begin{table}[t!] 
\TABLE
{Summary of different approaches to generate $\{x_t,u_t\}_{t=1}^T$.\label{tab:rate-comparison}}
{\begin{tabular}{l | c | c | c }
Approach & Binary search & No. iterations & Per-iteration cost\\ \hline
OFO-based & $\log(1/\epsilon)$ & $\max\left\{ r_u(\epsilon), r_x(\epsilon) \right\}$ & $O(km + mn)$\\
FO-based pessimization & $\log(1/\epsilon)$ & $r_x(\epsilon)$ & $m \Pess(\epsilon,k) + O\left(mn\right)$\\
Nominal oracle & see Table \ref{tab:nom-complexity} & $r_u(\epsilon)$ & $O\left(km\right) + \Nom(\epsilon,m,n)$\\ \hline
Full pessimization & see Table \ref{tab:nom-complexity} & $O(1/\epsilon^n)$ & $m \Pess(\epsilon,k) + \Nom(\epsilon,m+t,n)^*$\\
& \multicolumn{3}{r}{\footnotesize $^*$(number of constraints is $m+t$ as it grows by at least 1 each iteration $t$)}\\
Direct FOM via CoMirror & $1$ & $O(1/\epsilon^2)$ & $m \Pess(\epsilon,k) + O\left(mn\right)$
\end{tabular}}
{}
\end{table}

\begin{table}[ht]
\TABLE
{Arithmetic complexity for different implementations of pessimization oracles.
\label{tab:pess-complexity}}
{
\begin{tabular}{l|c}
Implementation & $\Pess(\epsilon,k)$\\ \hline
Closed form & $O(k)$ \\
IPM & $O(k^3\log(1/\epsilon))$ \\
FOM$^*$ & $O\left(k\log(1/\epsilon)\right)$ \\
 \multicolumn{2}{c}{\footnotesize$^*$(when $f^i$ are smooth, strongly convex in $u^i$)}
\end{tabular}
}
{}
\end{table}

\begin{table}[ht]
\TABLE
{Arithmetic complexity for different implementations of nominal oracles.\label{tab:nom-complexity}}
{
\begin{tabular}{l|ccc}
Implementation & $\Nom(\epsilon,m,n)$ & Type & Binary search\\ \hline
IPM & $O\left(km + \sqrt{m} (n^3 + mn) \log(1/\epsilon) \right)$ & optimization & 1\\
CoMirror & $O\left(mn/\epsilon^2\right)$ & optimization & 1\\
Convex-concave SP$^*$ & $O\left(\log(m)mn/\sqrt{\epsilon} \right)$ & feasibility & $\log(1/\epsilon)$\\
 \multicolumn{2}{l}{\footnotesize \quad$^*$(when $f^i$ are smooth, strongly convex in $x$)} & &
\end{tabular}
}
{}
\end{table}

Recall that Table~\ref{tab:rate-comparison} summarizes the rates for the various approaches, together with rates for the full pessimization approach of \cite{MutapcicBoyd2009} and using the CoMirror with pessimization (discussed in Section \ref{sec:connections-other-FOM}). Note that the total \emph{overall arithmetic complexity} of each approach is obtained by multiplying the quantities in each row in Table~\ref{tab:rate-comparison}. The quantities $r_u(\epsilon), r_x(\epsilon)$ will generally be $O(1/\epsilon^2)$, with potential for application-specific acceleration when the functions $f^i$ exhibit favorable structure. 
Table~\ref{tab:rate-comparison} indicates that our FO-based pessimization approach when it admits a closed form solution for the implementation of $\Pess(\epsilon,k)$ and the nominal feasibility oracle-based approach which uses a polynomial-time IPM solver to implement the nominal feasibility oracle $\Nom(\epsilon,m,n)$ give the best dependence on $\epsilon$ among all of the methods. These are better than our OFO-based approach by factors of $\max\{1,r_u(\epsilon)/r_x(\epsilon)\}$ and $\max\{1,r_x(\epsilon)/r_u(\epsilon)\}$ respectively. However, in many applications, we can expect that $r_u(\epsilon) \approx r_x(\epsilon)$, so these factors will be constant. In this case, our OFO-based approach becomes competitive with having a closed form pessimization oracle in our FO-based pessimization approach or using a nominal IPM solver in \cite{BenTalHazan2015}. That said, compared to IPMs, our OFO-based approach demands much less memory, and it is able to maintain a much lower dependence on the dimensions $m,n$ and thus is much more scalable, whereas the cost per iteration of such IPMs has a rather high dependence on the dimension. In addition, the memory requirements of IPMs are far more than OFO algorithms, posing a critical disadvantage to their use in large-scale applications. Similar comparisons of our OFO-based approach against pessimization or nominal feasibility oracle-based approaches utilizing other methods point out its advantage, which is at least an order of magnitude better in terms of its dependence on $\epsilon$. In fact, when $r_x(\epsilon)\approx r_u(\epsilon)$, our method can lead to savings over the approach of \cite{BenTalHazan2015} with CoMirror algorithm used in its oracle by a factor as large as $O(1/(\epsilon^2 \log(1/\epsilon)))$.

\subsection{Connections with Existing First-order Methods}\label{sec:connections-other-FOM}

Finally, we would like to discuss and contrast directly solving robust convex optimization problems \eqref{eqn:robust-opt-intro} via general first-order methods. 
Many FOMs require domains that are simple so that the prox operations can be easily done. In that respect, domains defined by multiple functional constraints $g^i(x) \leq 0$ creates a challenge for directly applying many of these algorithms. We now discuss two existing classes of FOMs that are designed to handle such domains: primal-dual methods and the CoMirror approach. Applying these FOMs to the RO problem \eqref{eqn:robust-opt-intro} can be viewed as another alternative solution methodology to solve RO problems without using the robust counterpart.

A general technique to address the functional constraints in the domain is to embed these constraints into the objective through Lagrange multipliers, and then solve the associated dual problem via FOMs (see e.g., \cite{NedicOzdaglar2009a}). Such methods are known as primal-dual methods. For the RO problem \eqref{eqn:robust-opt-intro}, this corresponds to solving
\[ \max_{\lambda} \left\{ \CL^*(\lambda) := \min_{x \in X} \left[ f^0(x) + \sum_{i=1}^m \lambda^{(i)} g^i(x) \right] : \lambda \geq 0 \right\}, \]
where we define $g^i(x) := \sup_{u^i \in U^i} f^i(x,u^i)$. Primal-dual methods (e.g., \cite{NedicOzdaglar2009a}) commonly require us to solve the inner minimization problem over $x \in X$ at each iteration. For RO, this means we must solve an expensive SP problem at each iteration. Our OFO-based approach aims to improve on this by reducing the per-iteration cost of each step to simple first-order updates. Two exceptions within the primal-dual methods are the work of Nedi\'{c} and Ozdaglar \cite{NedicOzdaglar2009b} and Yu and Neely \cite{YuNeely2016}, which have cheap per-iteration cost based on only gradient computations and projection operations in the Euclidean setup. Nedi\'{c} and Ozdaglar \cite{NedicOzdaglar2009b} provide a convergence rate of $O(1/\sqrt{T})$ in the non-smooth case. While using such a primal-dual method has the advantage that no binary search is needed, we note that this requires two assumptions to guarantee convergence: we have access to exact first-order information for the robust constraint functions $g^i(x) := \sup_{u^i \in U^i} f^i(x,u^i)$, and the standard Slater constraint qualification condition (i.e., strict feasibility) is satisfied. The first assumption is often not satisfied, since we may only be able to compute $g^i(x)$ up to accuracy $\epsilon$. While there exists some FOMs that work with inexact objective gradients over simple domains, see e.g., \cite{DevolderNesterovGlineur2014}, such methods have only been applied to specific max-type objectives, e.g., objectives obtained from smoothing. It is unclear how such methods can be extended for more general max-type functions which can arise in RO. Secondly, enforcing the Slater condition implicitly enforces feasibility of \eqref{eqn:robust-opt-intro}. In contrast, our framework directly uses the functions $f^i(x,u^i)$, so it does not need to take into account the inexact gradient information, and can certify infeasibility of \eqref{eqn:robust-opt-intro}.
Yu and Neely \cite{YuNeely2016} present a method that can guarantee $O(1/T)$ convergence when all functions are smooth. However, for RO problems, the constraint functions $g^i(x)$ are non-smooth due to the supremum operation, thus their results do not apply to RO.

The only FOM that we are aware of that can solve convex problems with functional constraints without assuming feasibility is the CoMirror algorithm \cite{BeckBenTal2010} and its earlier variations in the Euclidean setup \cite{NemYudin1983,Nesterov2004Book,Polyak1967}. The CoMirror\footnotemark[1] algorithm finds an $\epsilon$-optimal $\epsilon$-feasible solution in $O(1/\epsilon^2)$ iterations to a convex program $\min_{x \in X} \left\{ f^0(x) : g^i(x) \leq 0,\ i \in [m] \right\}$ or certifies its infeasibility by using (sub)gradient information of the objective $f^0$ as well as the constraint functions $g^i$.
\footnotetext[1]{Recall that the CoMirror algorithm is also discussed in Section~\ref{sec:RateDiscussion} as a method to implement the nominal feasibility solver; in that case we are given the noises $\bar{u}^i$ resulting in $g^i(x) := f^i(x,\bar{u}^i)$, and thus the subgradient of $g^i(x)$ is simply the subgradient of $f^i(x,\bar{u}^i)$.}
In the RO problem \eqref{eqn:robust-opt-intro} we defined $g^i(x) := \sup_{u^i \in U^i} f^i(x,u^i)$. As mentioned above, in many cases, we may only be able to compute $g^i(x)$ approximately, thus only have access to approximate/inexact gradient information. It is unknown to us whether or not techniques such as the ones from \cite{DevolderNesterovGlineur2014} can be applied to the CoMirror algorithm in the presence of this type of gradient information. While the CoMirror algorithm's complexity is $O(1/\epsilon^2)$  (see also Nesterov \cite[Chapter 3.2.4]{Nesterov2004Book} for a similar result in the Euclidean case), our iterative framework can exploit favorable structure on the functions $f^i$ that can improve on the iteration complexity $r_u(\epsilon), r_x(\epsilon)$. For the Euclidean case, Nesterov \cite[Chapters 2.3.4-2.3.5]{Nesterov2004Book} shows also that convergence can be obtained in $O(\log(1/\epsilon))$ iterations when the objective and all constraint functions are both smooth and strongly convex in $x$. However, such an improvement does not apply to the RO problem, since we cannot in general guarantee that $g^i(x)=\sup_{u^i \in U^i} f^i(x,u^i)$ is smooth in $x$. It is unknown whether the iteration complexity of CoMirror algorithm can be improved when only the underlying function $f^i(x,u^i)$ is strongly convex or smooth, or when $g^i(x)$ is strongly convex but non-smooth.

Finally, let us get back to the case when we have a robust feasibility problem with a \emph{single} constraint $m=1$ and a convex uncertainty set $U=U^1$. In such a case, as discussed in Remark~\ref{rem:SPm=1}, we have a direct convex-concave SP problem \eqref{eqn:SadVal} under Assumption~\ref{ass:f-concave-u}. The OFO-based approach then corresponds to bounding the whole SP gap \eqref{eqn:SPgap}, the FO-based pessimization corresponds to bounding the primal gap i.e., the first term in \eqref{eqn:SPgap}, and the nominal feasibility oracle approach corresponds to bounding the dual gap, i.e., the second term in \eqref{eqn:SPgap}. Without any further structural assumptions on $f^1$, convex-concave SP problems can be solved in $O(1/\epsilon^2)$ iterations. Our approaches also achieve this rate immediately, see  Table~\ref{tab:rate-comparison}. Moreover, when $f^1$ is smooth in $x$ and strongly concave in $u^1$, our general framework can achieve a rate of $O(1/\epsilon)$.  
However, for specific applications involving a single robust constraint, directly working with a specialized convex-concave SP formulation can improve this rate further. For example, such an improved rate of $O(1/\sqrt{\epsilon})$ is achieved in \cite{Ben-TalRobustSVM2012} for a robust support vector machine problem.

\section{Application Example: Robust Quadratic Programming}\label{sec:Application}
Our framework is general and can be applied to many robust convex optimization problems. In this section we walk through the setup and resulting convergence rates of our framework for a robust feasibility problem of a quadratically constrained quadratic program (QP) with ellipsoidal uncertainty. To be precise, our deterministic feasibility problem is 
\begin{align*}
\text{find}~~ x\in X~~\text{s.t.}~~ \|A_i x\|_2^2 \leq b_i^\top x +c_i,\quad \forall i\in[m],
\end{align*}
where $X \subseteq \R^n$ is the unit Euclidean ball, $A_i\in\R^{n\times n}$, $b_i\in\R^n$, and $c_i\in\R$ for all $i\in[m]$. 
We consider the robust quadratic feasibility problem given by
\begin{equation}\label{eqn:robustQP-feas} 
\text{find}~~ x\in X~~\text{s.t.}~~\sup_{u \in \widehat{U}} \left\| \left( A_i + \sum_{k=1}^{K} P_k^i\, u^{(k)} \right) x \right\|_2^2 - b_i^\top x -c_i \leq 0,\quad \forall i\in[m], 
\end{equation}
where $P_1^i,\ldots,P_K^i$ are uncertainty matrices for each constraint $i\in[m]$, for simplicity we assume uncertainty sets $U^i = \widehat{U}= \{u \in \R^{K} : \|u\|_2 \leq 1\}$ for all $i\in[m]$, and $u^{(k)}$ denotes the $k$-th entry of $u$. 

It is well known that the robust counterpart of this feasibility problem is a semidefinite program \cite{BenTalelGhaouiNemirovski2009,BertsimasBrownCaramanis2011}. Because current state-of-the-art QP solvers can handle two to three orders
of magnitude larger QPs than semidefinite programs (SDPs), Ben-Tal et al.\@ \cite[Section 4.2]{BenTalHazan2015} suggest an approach that avoids solving SDPs associated with robust QPs. Their approach relies on running a probabilistic OCO algorithm in which a trust region subproblem (TRS)---a class of well-studied nonconvex QPs---is solved in each iteration. 
Our results here further enhance this approach. In particular, we show that we can achieve the same rate of convergence in our framework while working with a deterministic OCO algorithm and only carrying out first-order updates in each iteration. In fact, the most expensive operation involved with each iteration of our approach is a maximum eigenvalue computation. Because maximum eigenvalue computation is much cheaper than solving a TRS, we not only present a deterministic approach but also strikingly reduce the cost of each iteration.

To simplify our exposition, let us introduce some notation. For each $i\in[m]$, we define the matrix $\cP_x^i  \in \R^{n \times K}$ whose columns are given by the vectors $P_k^i \, x$ for $k \in [K]$ together with 
\begin{align*}
Q_x^i := (\cP_x^i)^\top \cP_x^i \in \Se^K_+,\quad
r_x^i := (\cP_x^i)^\top A_i x \in \R^K, \quad \text{and}\quad
s_x^i := \|A_i x\|_2^2 - b_i^\top x -c_i \in \R;
\end{align*}
then it is easy to check that for all $i\in[m]$ and $u\in\R^K$ we have
\[ \left\| \left( A_i + \sum_{k=1}^{K} P_k^i\, u^{(k)} \right) x \right\|_2^2 - b_i^\top x -c_i = u^\top Q_x^i u + 2(r_x^i)^\top u + s_x^i. \]
For each $i\in[m]$, we define $f^i: X \times \widehat{U} \to \R$ as
\begin{align}
f^i(x,u) &:= \left\| \left( A_i + \sum_{k=1}^{K} P_k^i\, u^{(k)} \right) x \right\|_2^2  - b_i^\top x -c_i + \lambda_{\max}(Q_x^i)\left( 1 - \|u\|_2^2 \right) \nonumber\\
&= u^\top Q_x^i  u + 2(r_x^i)^\top u + s_x^i + \lambda_{\max}(Q_x^i)\left( 1 - \|u\|_2^2  \right). \label{eqn:f-forQP}
\end{align}

\begin{lemma}\label{lem:QPreformulation}
For each $i\in[m]$, the function $f^i(x,u)$ defined in \eqref{eqn:f-forQP} is convex in $x$ for any fixed $u\in  \widehat{U}$ and concave in $u$ for any given $x$. Moreover, for all $i\in[m]$ and for any $x\in X$,
\[ \sup_{u \in \widehat{U}} \left\| \left( A_i + \sum_{k=1}^{K} P_k^i\, u^{(k)} \right) x \right\|_2^2 - b_i^\top x -c_i =\sup_{u \in \widehat{U}} f^i(x,u) . \]
\end{lemma}
\proof{Proof.}
Fix $i\in[m]$. 
By rearranging terms in \eqref{eqn:f-forQP}, we obtain $f^i(x,u)=u^\top (Q_x^i-\lambda_{\max}(Q_x^i) I_K) u + 2(r_x^i)^\top u +s_x^i$. Since $Q_x^i-\lambda_{\max}(Q_x^i) I_K\in\Se^K_+$ for any given $x$, $f^i(x,u)$ is concave in $u$ for any given $x$.

Now consider a fixed $u\in \widehat{U}$. Note that  
\[\lambda_{\max}(Q_x^i) = \max_{\|v\|_2 \leq 1} v^\top (Q_x^i) v = \max_{\|v\|_2 \leq 1} \sum_{1 \leq j,k \leq K} v^{(j)} v^{(k)} x^\top \!(P_j^i)^\top P_k^i\, x = \max_{\|v\|_2 \leq 1} x^\top\! \left(\! \sum_{k=1}^{K} P_k^i\, v^{(k)} \!\!\right)^\top\!\!\! \left(\! \sum_{k=1}^{K} P_k^i\, v^{(k)} \!\!\right) x.\]
Because $ \left( \sum_{k=1}^{K} P_k^i\, v^{(k)} \right)^\top \! \left( \sum_{k=1}^{K} P_k^i\, v^{(k)} \right)\in\Se^n_+$, then $\lambda_{\max}(Q_x^i)$ is a maximum of convex quadratic functions of $x$ and hence is convex in $x$. Thus, for fixed $u \in  \widehat{U}$, $f^i(x,u)$ is convex in $x$. 

Reformulation of the nonconvex QP over an ellipsoid into a convex QP over the ellipsoid 
via the relation between $u^\top\! Q_x^i  u + 2(r_x^i)^\top\! u + s_x^i$ and $f^i(x,u)$ in \eqref{eqn:f-forQP} follows from \cite[Theorem 2.1]{JeyakumarLi2013}.
\Halmos
\endproof

Lemma~\ref{lem:QPreformulation} implies that $\sup_{u \in \widehat{U}} f^i(x,u) \leq 0$ is an alternate representation of our robust quadratic constraint. We next state the convergence rate in our framework for the associated feasibility problem. For this, we define the quantities
\begin{align}\label{eqn:QPquantities}
&\sigma^2 := \max_{i\in[m]} \sum_{k=1}^{K} \|P_k^i\|_{\Fro}^2, \quad
\chi := \max_{i\in[m]} \max_{k \in [K]} \|P_k^i\|_{\Spec}, \quad \text{and}\quad \notag\\
&
\rho := \max_{i\in[m]}\|A_i\|_{\Spec}, \qquad~~
\beta := \max_{i\in[m]}\|b_i\|_2.
\end{align}
Note that $\chi\leq \sigma$. Furthermore, \cite[Lemma 7]{BenTalHazan2015} proves that $\|Q_x^i\|_{\Fro}\leq \sigma^2$ and  $\|r_x^i\|_2\leq \sigma\rho$ holds for all $x$ such that $\|x\|_2\leq 1$. 

\begin{corollary}\label{cor:robustQP}
Let our domain be given by $X = \{x\in\R^n:\; \|x\|_2 \leq 1\}$. 
The customization of our OFO-based approach to the problem~\eqref{eqn:robustQP-feas} ensures that within $O\left( ((\rho + \sqrt{K} \sigma)^2 + \beta)^2 \right) \epsilon^{-2}$ iterations, we obtain a robust feasibility/infeasibility certificate. 
Moreover, each iteration in our framework relies on a first-order update where the most expensive operation in the case of \eqref{eqn:robustQP-feas} is computing $\lambda_{\max}(Q_x^i)$, which can be done efficiently.
\end{corollary}
\proof{Proof.}
In order to apply OFO-based approach, we need to customize our proximal setup. Given that the sets $X$ and $ \widehat{U}$ are Euclidean balls, we set the proximal setup for generating the iterates $\{x_t,u^i_t\}_{t=1}^T$ to be the standard Euclidean \dgf\ with $\|\cdot\|_2$-norm, and thus $\Omega_X=\Omega_{\widehat{U}}={1\over 2}$. We must bound the magnitude of the gradients measured by the $\|\cdot\|_2$-norm.  
Note that for any $i\in[m]$, the gradients of $f^i$ are given by
\begin{align*}
\grad_{u} f^i(x,u) &= 2 \left( Q_x^i - \lambda_{\max}(Q_x^i) I_K \right) u + 2 r_x^i\\
\grad_x f^i(x,u) &= 2 \left(\!\! A_i + \sum_{k=1}^{K} P_k^i\, u^{(k)} \!\!\right)^\top \!\!\left(\!\! A_i + \sum_{k=1}^{K} P_k^i\, u^{(k)} \!\!\right) x + 2 \left(1 - \|u\|_2^2\right)\! \left(\! \sum_{k=1}^{K} P_k^i\, v^{(k)} \!\right)^\top\!\! \left(\! \sum_{k=1}^{K} P_k^i\, v^{(k)} \!\right) x -b_i,
\end{align*}
where $v \in  \widehat{U}$ is an eigenvector of $Q_x^i$ corresponding to $\lambda_{\max}(Q_x^i)$.

Let us fix an $i\in[m]$. We first bound $\|\grad_u f^i(x,u)\|_2$ for any $u\in\widehat{U}$ as follows:
\begin{align*}
\|\grad_u f^i(x,u)\|_2 &= 2 \left\|  \left( Q_x^i - \lambda_{\max}(Q_x^i) I_K \right) u +  r_x^i \right\|_2\\
&\leq 2\left(\left\| \left( Q_x^i - \lambda_{\max}(Q_x^i) I_K \right) u\right\|_2 + \left\| r_x^i \right\|_2\right)
\leq 2 \lambda_{\max}(Q_x^i) \|u\|_2 + 2 \sigma \rho
\leq 2 (\sigma^2 +  \sigma \rho),
\end{align*}
where the second inequality follows from $\left\| Q_x^i - \lambda_{\max}(Q_x^i) I_K \right\|_{\Spec} \leq \lambda_{\max}(Q_x^i)$ and $\| r_x^i\|_2\leq \sigma \rho$ which is implied by \cite[Lemma 7]{BenTalHazan2015},  and the last inequality follows from the facts that $u\in\widehat{U}$, the definitions given in \eqref{eqn:QPquantities}, and $\lambda_{\max}(Q_x^i) = \|\cP_x^i\|_{\Spec}^2 \leq  \|\cP_x^i\|_{\Fro}^2 \leq \sum_{k=1}^{K} \|P_k^i\|_{\Fro}^2 \leq \sigma^2$ for any $x\in X$. 
Therefore, we deduce from Theorem~\ref{thm:OCO-non-smooth} with uniform weights $\theta_t= 1/T$ that the rate of convergence for bounding the weighted regret associated with constraint $i\in[m]$ using the online mirror descent algorithm is 
\[ \sup_{u \in U} \frac{1}{T} \sum_{t=1}^{T} f^i(x_t,u) - \frac{1}{T} \sum_{t=1}^{T} f^i(x_t,u_t) \leq \frac{2(\sigma^2 + \sigma \rho)}{\sqrt{T}}.\]
This implies that $r_u(\epsilon) = O((\sigma^2 + \sigma \rho)^2 \epsilon^{-2})$.

We next bound the weighted regret of the functions $\varphi_t(x) = \max_{i \in [m]} f^i(x,u_t^i)$, i.e., the term $\epsilon^\bullet(\{x_t,u_t,\theta_t\}_{t=1}^T)$ by bounding the $\|\cdot\|_2$-norm of $\grad_x \varphi_t(x)$. Notice that
\begin{align*}
\left\| \grad_x \varphi_t(x) \right\|_2 &\leq \max_{i \in [m]} \|\grad_x f^i(x,u_t^i)\|_2.
\end{align*}
Thus, we must bound $\|\grad_x f^i(x,u)\|_2$ for all $x \in X$, $u \in \widehat{U}$. To this end, note that for any $u\in\widehat{U}$
\[ \left\| \sum_{k=1}^{K} P_k^i\, u^{(k)} \right\|_{\Spec} \leq \sum_{k=1}^{K} \|P_k^i\|_{\Spec}\, |u^{(k)}| \leq \sqrt{K} \max_{k \in [K]} \|P_k^i\|_{\Spec} \leq \sqrt{K} \chi, \]
where the second inequality holds because $\|u\|_1\leq\sqrt{K}\|u\|_2\leq\sqrt{K}$ holds for all $u\in\widehat{U}$. 
Then for any $x\in X$, $u\in\widehat{U}$, and eigenvector $v\in\widehat{U}$, we have
\begin{align*}
\|\grad_x f^i(x,u)\|_2 
&\leq 2\left\| A_i + \sum_{k=1}^{K} P_k^i\, u^{(k)}\right\|_{\Spec}^2 \|x\|_2 + 2\left(1 - \|u\|_2^2\right)\left\| \sum_{k=1}^{K} P_k^i\, v^{(k)} \right\|_{\Spec}^2 \|x\|_2 +\|b_i\|_2\\
&\leq 2 (\rho + \sqrt{K} \chi )^2 + 2 K \chi^2 +\beta\\
&\leq 4 (\rho + \sqrt{K} \sigma )^2 +\beta.
\end{align*}
Hence, $\left\| \grad_x \varphi_t(x) \right\|_2 \leq 4 (\rho + \sqrt{K} \sigma )^2 +\beta$. Then Theorem~\ref{thm:OCO-non-smooth} with weights $\theta_t=1/T$ implies 
\[ \sum_{t=1}^T \theta_t \max_{i \in [m]} f^i(x_t,u_t^i) - \inf_{x \in X} \sum_{t=1}^T \max_{i \in [m]} f^i(x,u_t^i) \leq \frac{\left( 4 (\rho + \sqrt{K} \sigma )^2 +\beta \right)}{\sqrt{T}}. \]
Thus, $r_x(\epsilon) = O\left( ((\rho + \sqrt{K} \sigma)^2 + \beta)^2 \right) \epsilon^{-2}$. Therefore, the number of iterations required for our OFO-based approach to obtain a robust feasibility/infeasibility certificate is $T = \max\{r_x(\epsilon),r_u(\epsilon)\}$.

Note that each iteration of our approach requires a first-order update that is composed of computing the gradients $\grad_x f^i(x,u)$ and $\grad_u f^i(x,u)$ and prox computations. Because our domains involve only direct products of Euclidean balls and simplices, they admit efficient prox computations which take $O(Km+mn)$ time. 
In order to evaluate the gradients $\grad_x f^i(x,u)$ and $\grad_u f^i(x,u)$, in addition to the elementary matrix vector operations, we need to compute $\lambda_{\max}(Q_x^i)$ which is the most expensive operation in our first-order update. Fortunately, computing the maximum eigenvalue of a matrix is a well-studied problem and can be computed very efficiently.
\Halmos
\endproof

In the case of robust QP feasibility problem~\eqref{eqn:robustQP-feas},  \cite[Corollary 3]{BenTalHazan2015} states that with probability $1-\delta$, their framework returns robust feasibility/infeasibility certificates in at most  \sloppy $O\left( K^2 \sigma^2(\rho^2 + \sigma^2) \log({m/\delta})\epsilon^{-2} \right)$ calls (iterations) to their oracle. In each call to their oracle, a nominal feasibility problem is solved to the accuracy ${\epsilon/2}$. In comparison we deduce from Corollary~\ref{cor:robustQP}  that our framework requires comparable number of iterations as the approach of Ben-Tal et al.\@ \cite{BenTalHazan2015}. Even so, there are a number of reasons that considerably favor our approach. First, our approach is deterministic as opposed to the high $1-\delta$ probability guarantee of \cite{BenTalHazan2015} which requires using an adaptation of the follow-the-perturbed-leader type OCO. Second, each iteration of their approach requires solving a nominal feasibility problem for solution oracle as well as solving TRSs for the computation of noises $u_t$. In contrast to this, in each iteration we carry out mainly elementary operations such as matrix vector multiplications and our most computationally expensive operation is the maximum eigenvalue computations $\lambda_{\max}(Q_x^i)$. While there are established algorithms to solve the TRS, it is inherently more complicated than finding the maximum eigenvalue of a positive semidefinite matrix. Moreover, \cite{BenTalHazan2015} suffers from the additional computational cost of their solution oracle which solves the nominal feasibility problem. Hence, our approach, while requiring a comparable number of iterations, reduces the cost per iteration remarkably.

\section{Numerical Study}\label{sec:numerical}

In this section, we conduct a numerical study comparing the approaches discussed so far. We consider the following quadratic program inspired by mean-variance portfolio optimization problems with a factor model for the return vector (see, e.g., \cite{GoldfarbIyengar2003}):
\begin{equation}\label{eqn:port-opt}
\min_{x} \left\{ \| V x\|_2^2 + x^\top D x - \lambda \mu^\top x :~
x \in \Delta_n
\right\} ,
\end{equation}
where $\mu\in\R^n$ is the expected return vector, the term $x^\top (V^\top V+D) x$ captures the risk associated with the portfolio via a factor model, and $\lambda \geq 0$ represents the trade-off between the expected return of the portfolio and the risk associated with the portfolio.

In the robust formulation of \eqref{eqn:port-opt}, we consider the case where the true parameters $\mu \in \bbR^n$ and $V \in \bbR^{m \times n}$ belong to uncertainty sets $\cM$ and $\cV$ of form
\[ 
\CM := \left\{ \mu :~ \mu_0-\gamma \leq \mu,~  \mu \leq \mu_0 + \gamma \right\}, 
\quad  \CV := \left\{ V = V_0 + \sum_{k=1}^K P_k u_k :~ \|u\|_2 \leq 1 \right\},
\]
where the nominal data $\mu_0 \in \R^n$, $\gamma\in\R^n$, and $V_0 \in \R^{m \times n}, \{P_k \in \R^{m \times n}\}_{k=1}^K$ are given to us. Then the robust problem is given by 
\begin{equation}\label{eqn:robust-port-opt}
\min_{x} \left\{ \max_{V \in \CV} \|Vx\|_2^2 + x^\top D x - \lambda \min_{\mu \in \CM} \mu^\top x :~
x \in \Delta_n
\right\}.
\end{equation}

Our test instances are synthetically generated, largely following the random instance generation model from \cite{GoldfarbIyengar2003}. We begin by specifying three parameters: $n$, the number of variables; $m$, the number of factors (which controls the rank of $V$); and $\alpha \in (0,1)$, a parameter controlling the size of the uncertainty sets. For each instance, we randomly generate matrices $V \in \bbR^{m \times n}$ and $F \in \bbR^{m \times m}$, where we ensure $F$ is positive semidefinite, and define $D = 0.1 \Diag(V^\top F V)$. We then generate $p > m$ factor samples $f_{(l)} \in \bbR^m$, $l \in [p]$, where each $f_{(l)} \sim N(0,F)$, and we also generate $\mu \in \bbR^n$ where each entry $\mu_i \sim U(1,5)$. We then set $\mu_{(l)} = \mu + V^\top f_{(l)} + \epsilon_l$, where $\epsilon_{(l)} \sim N(0,D)$ are independent of the factor sample $f_{(l)}$. The matrices $\mu$ and $V$ are estimated via linear regression on $\mu_{(l)}$ and $f_{(l)}$, to obtain $\bar{\mu}, \bar{V}$. The nominal data for \eqref{eqn:port-opt} are set to be $\mu_0 = \bar{\mu}$, $V_0 = F^{1/2} \bar{V}$. To define the uncertainty sets, we first compute the scaled sum of squared errors for each $i \in [n]$, $s_i^2 = \frac{1}{p-m-1} \sum_{l=1}^{p} (\mu_{(l),i} - \mu_{0,i} - V_{0,i}^\top f_{(l)})^2$. Let $c_J(\alpha)$ be the $\alpha$-critical value of an $F$-distribution with $J$ degrees of freedom, and let $\nu$ be the top-left entry of $A^{-1}$, where $A \in \bbR^{(p+1) \times (p+1)}$ is the Gram matrix of the vectors $\bm{1}_m, \{f_{(l)}\}_{l=1}^p$. Then we set $\gamma_i = \sqrt{\nu c_1(\alpha) s_i^2}$ for $i\in[n]$, which defines the uncertainty set for $\mu$. The uncertainty set for $V$ is chosen by randomly generating matrices $P_k$, and then scaling them appropriately so that the norm of each column $i$ of $V-V_0$ is at most $\sqrt{m c_m(\alpha) s_i^2}$ for every $V \in \CV$.

We set $p=90$ and $\alpha = 0.95$, while varying $m\in\{3,5,7,10,15,20,25\}$ and $n\in\{100,200,300,400,500,600,700\}$. We fix the underlying dimension of the uncertainty set $\cV$ to be  $K = \min\{2m,15\}$. We generate five instances for each combination of $m$ and $n$.

The four approaches we test are our OFO-based approach from Section \ref{sec:RegretforRO}, our FO-based pessimization approach from Section \ref{sec:pessimization-oracle} (see Theorem~\ref{thm:pessimization-oracle-approach}), the nominal oracle-based approach of \cite{BenTalHazan2015} from Section \ref{sec:oracle-based-approach}, and the full pessimization approach of \cite{MutapcicBoyd2009}, which requires both a pessimization and an extended nominal feasiblity oracle. Since \eqref{eqn:robust-port-opt} is an instance of a robust quadratic program, the form for nominal and pessimization oracles can be derived from Section \ref{sec:Application}. One-dimensional line search using Brent's algorithm \cite{Brent1973} was used to choose step sizes for each iteration of FO-based methods. An error tolerance of $\epsilon=0.002$ is used in all instances.

Experiments are performed on a Linux machine with 2.8GHz processor and 64GB memory using Python v3.5.2. Whenever the nominal (extended nominal) oracles and pessimization oracles do not have closed form solutions, they are implemented in Gurobi v7.0.2. We use standard Gurobi tolerances and parameter choices. We employ the implementation of Brent's algorithm in Python's \texttt{scipy.optimize} package.

\begin{figure}[h]
\centering
\includegraphics[width=100mm]{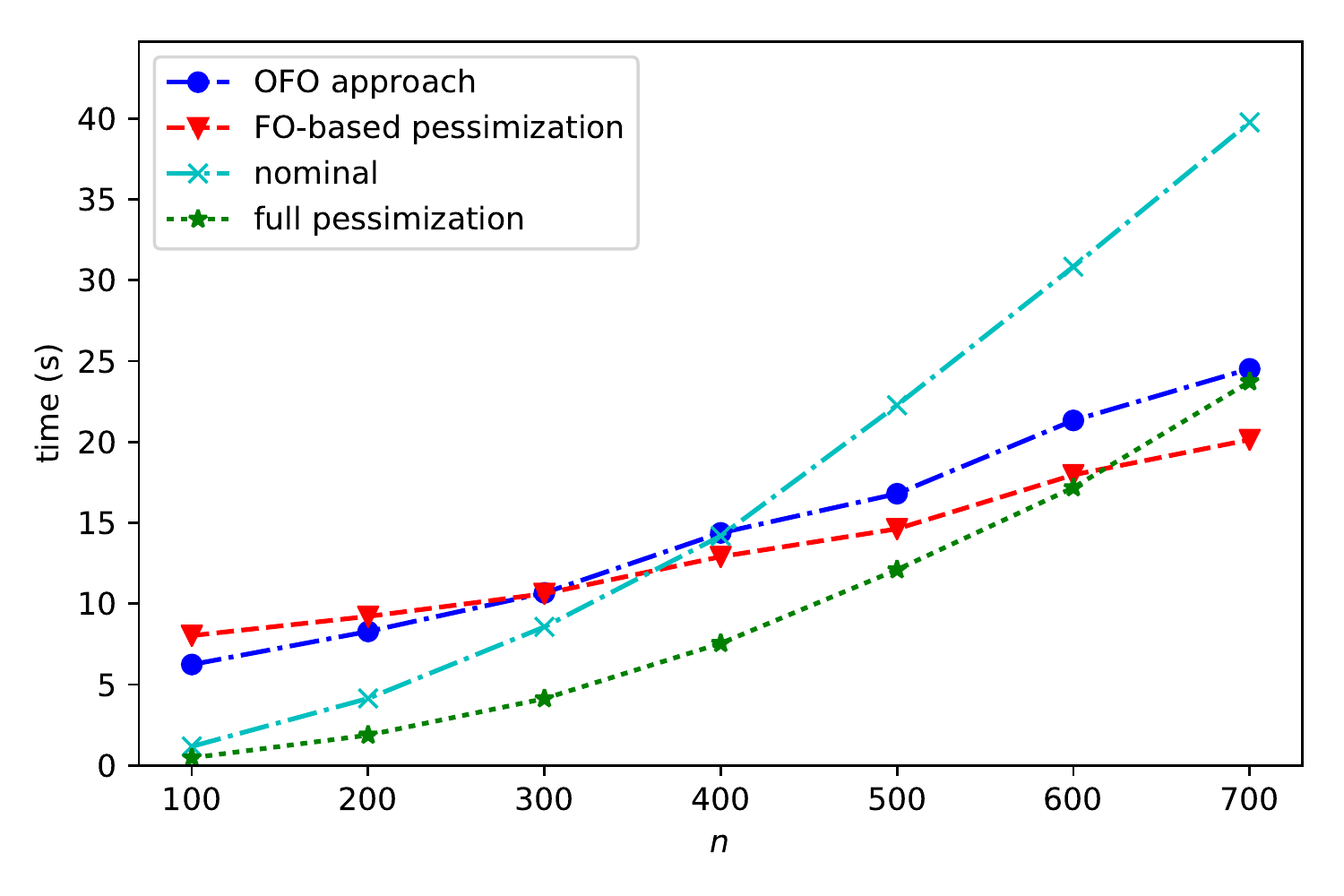}
\vspace{-10pt}
\caption{Average solve times (seconds) for different $n$.
}
\label{fig:average-solve-times}
\end{figure}

Figure \ref{fig:average-solve-times} plots the average solve times in seconds against different $n$ for each of the approaches, averaging across all $m$. As we expect, for low dimensions $n$, the oracle-based approaches solve the instances very quickly compared to our first-order based approaches. However, as $n$ increases to $400,500,600,700$, we see that the solution times of our first-order based approaches beat the nominal oracle approach, and become comparable to the full pessimization approach for $n=600,700$. In particular, we observe that the FO-based pessimization approach (see Theorem~\ref{thm:pessimization-oracle-approach}) solves faster when $n=700$.

\begin{figure}[h]
\centering
\subfigure[$n=400$]{
\includegraphics[width=78mm]{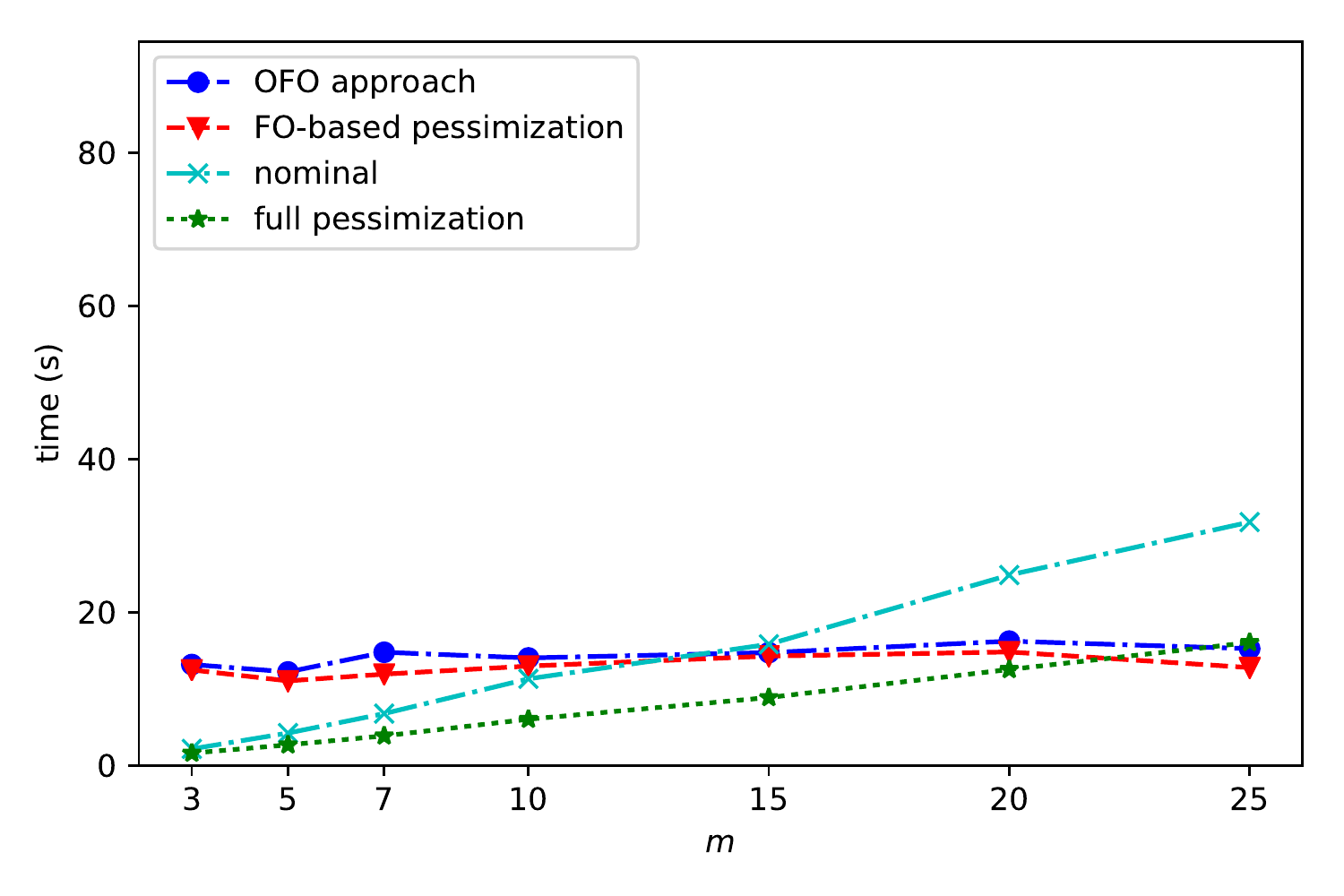}
}
\subfigure[$n=500$]{
\includegraphics[width=78mm]{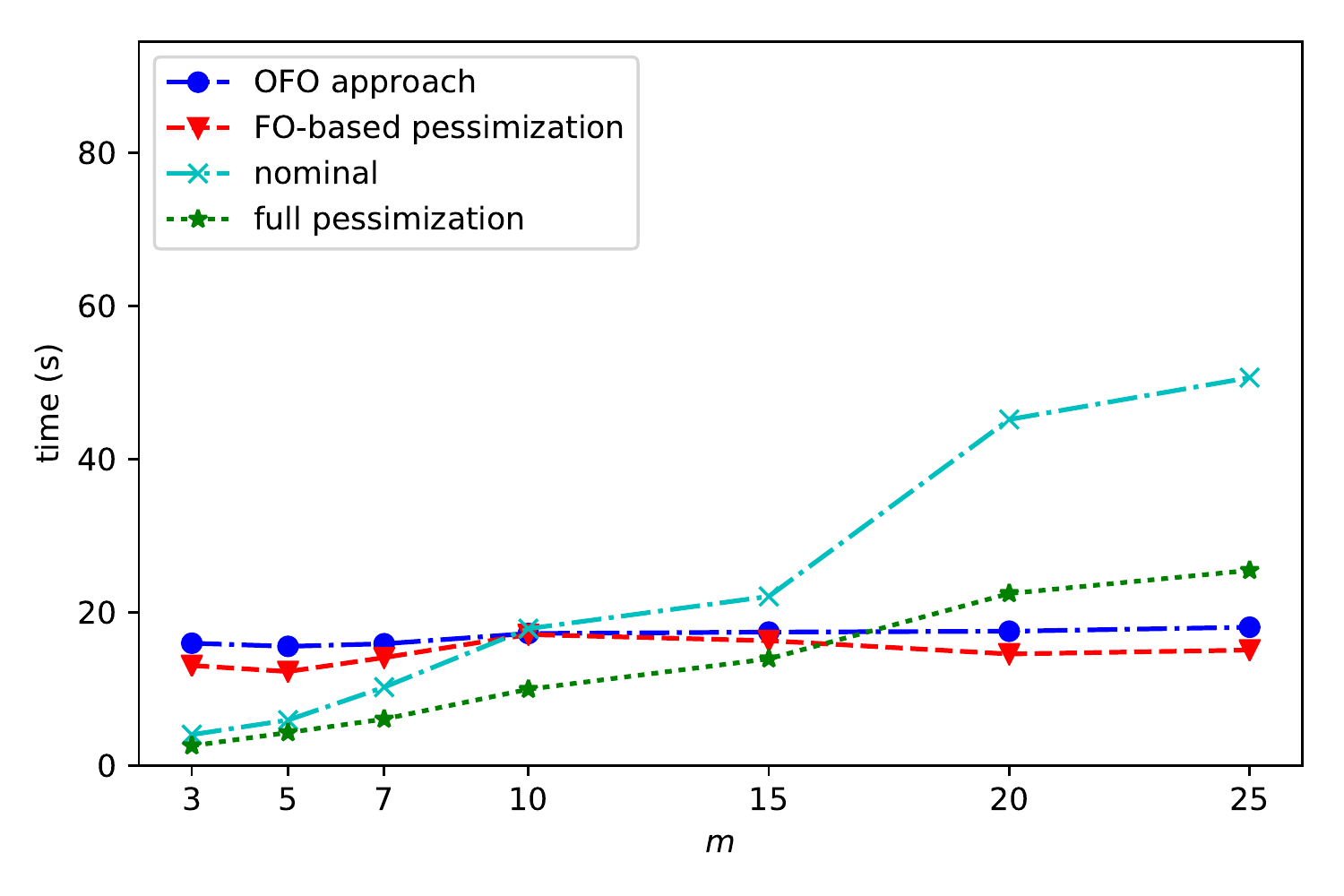}
}
\subfigure[$n=600$]{
\includegraphics[width=78mm]{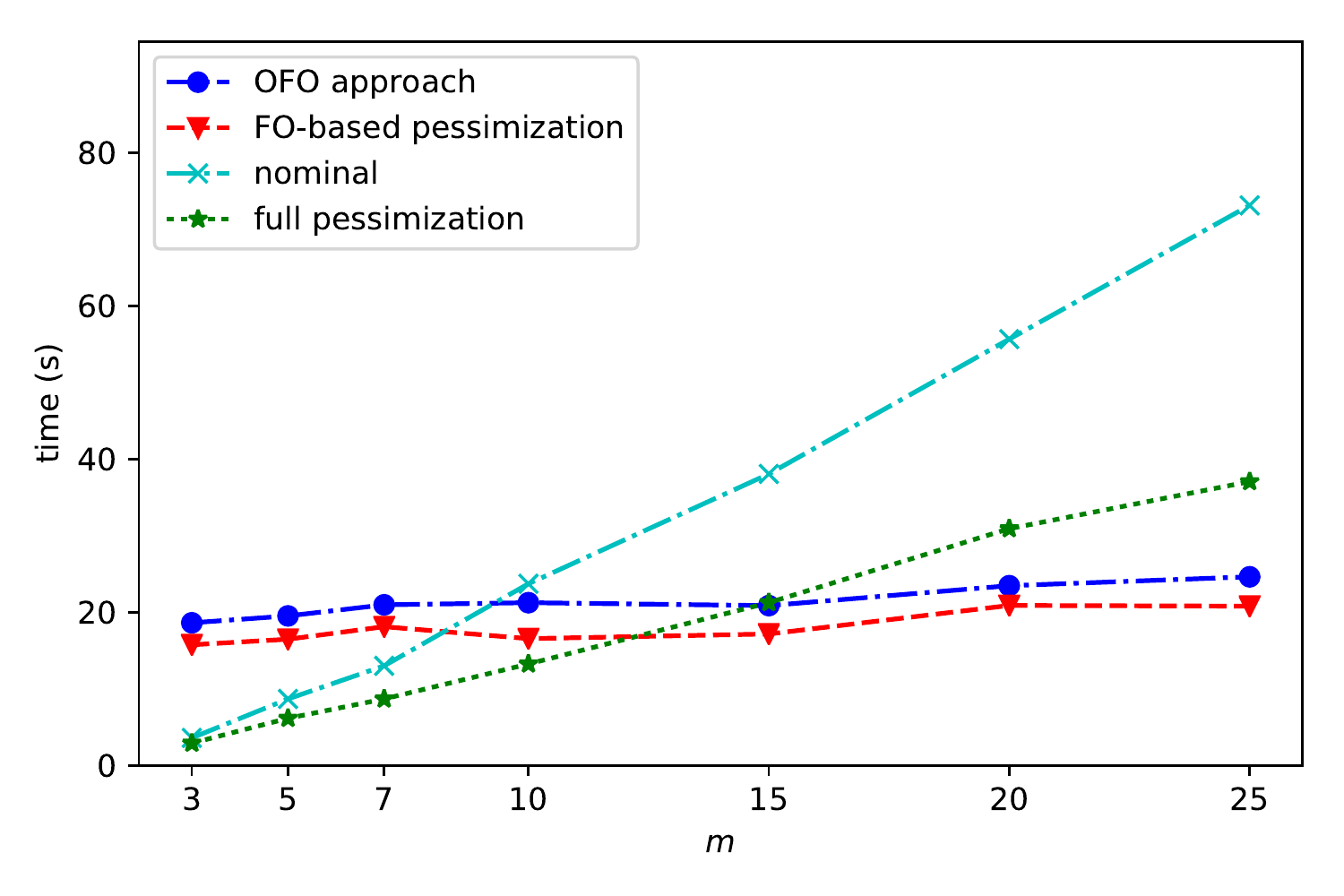}
}
\subfigure[$n=700$]{
\includegraphics[width=78mm]{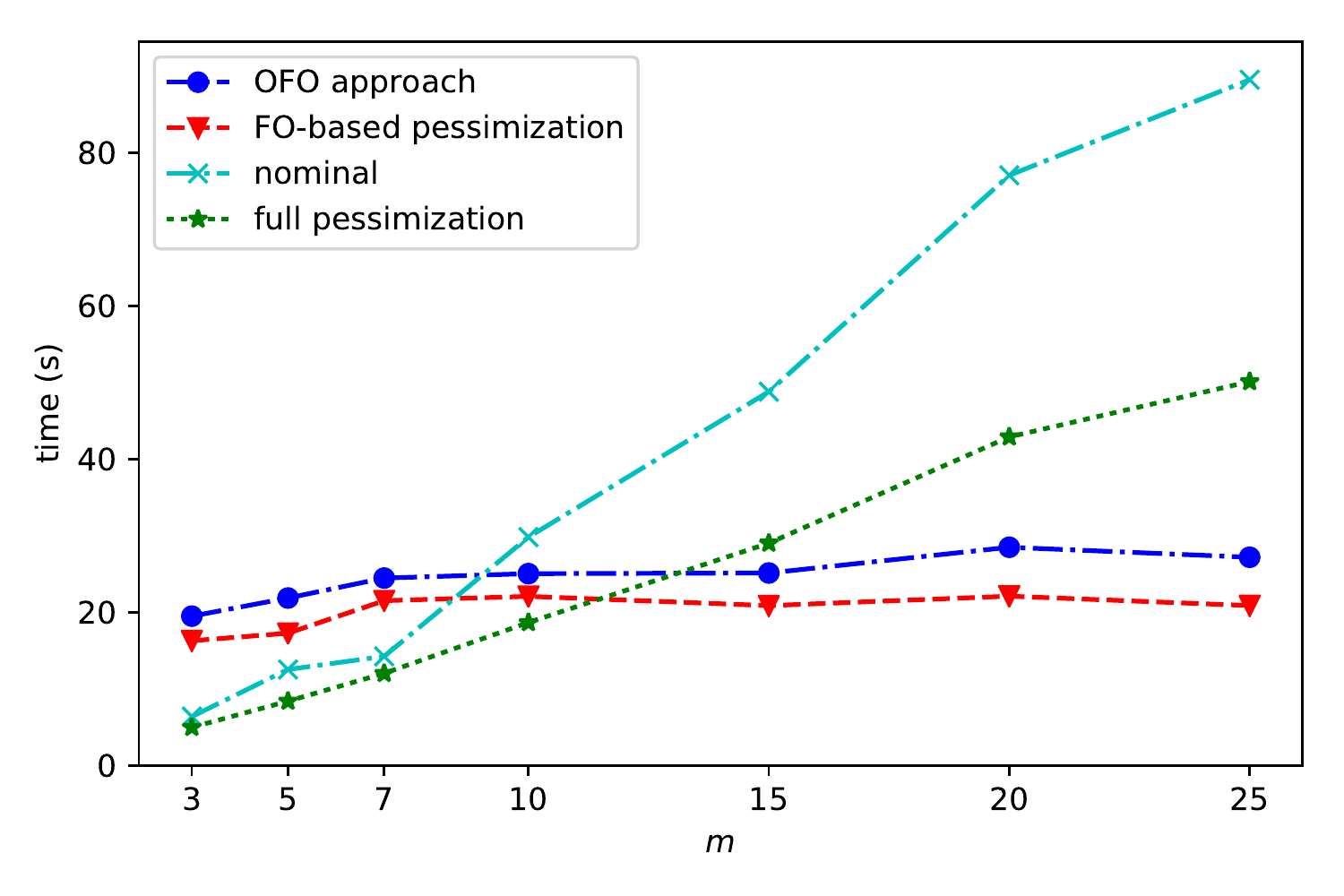}
}
\caption{Average solve times (seconds) for different $n$ and $m$.}
\label{fig:average-solve-times-plots}
\end{figure}

The dimension $m$ influences the rank of the nominal matrix $V_0^\top V_0$ and controls the difficulty of the problems. Examining the results for different $m$ further highlights the benefits of utilizing the first-order based approaches. Figure \ref{fig:average-solve-times-plots} plots average solve times for different $m$ while fixing $n = 400,500,600,700$. For the oracle-based methods, the solution times increase with $m$, while the solution times for first-order based methods remains relatively constant with $m$. For $m \geq 20$, we observe that our first-order based approaches significantly outperforms the oracle-based methods which require a nominal solver.
Notice that, while we expect our OFO-based approach to outperform the FO-based pessimization approach due to the burden of solving an eigenvalue problem in each iteration for computing the pessimization oracle, our results indicate the opposite. This is because for small values of $K$, calling a pessimization oracle is faster than the line search performed in the FO-based noise update. However, we believe that as $K$ increases, one-dimensional line search will become more efficient.

Finally, we examine the number of iterations and cost per iteration of different approaches averaged across all instances in Table \ref{tab:method-iters-summary}. We observe that, contrary to their theoretical iteration guarantees, the oracle-based approaches of \cite{MutapcicBoyd2009,BenTalHazan2015} need very few iterations to find a solution. However, as expected, the average time per iteration is significantly higher for these methods due to their reliance on full nominal optimization solvers. 
This further highlights the benefit of utilizing first-order methods for robust optimization when the deterministic version of the problem is already very expensive, and hence nominal oracles become expensive.
\begin{table}[h!tb]
\TABLE
{Average number of iterations and average time per iteration for each approach.
\label{tab:method-iters-summary}}
{\vspace{-10pt}
\begin{center}
\begin{tabular}{l|rc}
\toprule
		&     \# iterations &  seconds per iteration \\ \hline
first-order        &  961.487 &          0.015 \\
FO-based pessimization      & 1009.054 &          0.013 \\
nominal            &    3.708 &          4.841 \\
full pessimization &    1.975 &          4.875 \\
\bottomrule
\end{tabular}
\end{center}
}
{}
\end{table}

\section{Conclusion}\label{sec:Conclusions}
In this paper, we advance the line of research in \cite{BertsimasSim2003,MutapcicBoyd2009,BenTalHazan2015} that aims to solve robust optimization problems via iterative techniques, i.e., without transforming them into their equivalent robust counterparts. Thus far, the literature on iterative methods for RO has relied on more expensive nominal feasibility or pessimization oracles. However, in many applications of robust convex optimization, the original deterministic problem comes equipped with first-order oracles that provide gradient/subgradient information on the constraint functions $f^i$. In this paper, we present an efficient framework that can both work with cheap online first-order oracles and also capture the prior oracle-based approaches of \cite{MutapcicBoyd2009} and \cite{BenTalHazan2015}. 
We further show that working with these OFO oracles essentially does not increase the worst case theoretical bound on number of overall oracle calls, i.e., the worst case bound on number of main iterations of our approach is better than or comparable to the prior approaches. Moreover, when OFO oracles are utilized in our framework, the resulting overall arithmetic complexity including all of the basic operations in each iteration is remarkably cheaper than the prior approaches. The resulting framework is simple, easy-to-implement, flexible, and it can easily be customized to many applications. We demonstrate  our framework via an illustrative robust QP example, where the most expensive operation in each iteration of our framework is a maximum eigenvalue computation. We further illustrate this with a preliminary numerical study on robust portfolio optimization problem. 

Our framework is amenable to exploiting favorable structural properties of the functions $f^i$ such as strong concavity, smoothness, etc., through which better convergence rates can be achieved. For example, when $f^i$ are strongly concave in $u^i$, by exploiting this structural information and using a customization of the weighted regret online mirror descent for strongly convex functions, it is possible to achieve a better convergence rate of $O(1/\epsilon)$ in both our online first-order oracle setup and the nominal feasibility oracle framework of \cite{BenTalHazan2015}. This then partially resolves/refines an open question stated in \cite{BenTalHazan2015} for the lower bound on the number of iterations/calls needed in their nominal feasibility oracle based framework. However, it remains open whether $O(1/\epsilon^2)$ bound is tight when no further favorable structure is present in $f^i$ or the tightness of $O(1/\epsilon)$ in the favorable case. 

There are several other compelling avenues for future research. From a practical perspective, it is well-known, and also confirmed by our preliminary proof-of-concept computational experiments, that the computation of gradients/subgradients constitute a major bottleneck in the practical performance of FOMs. Thus, as a step to reduce the  efforts involved in such computations, possible incorporation of stochastic \cite{RobbinsMonro1951,NemJudLanShapiro2009} and/or randomized FOMs \cite{JKKN13,BenTalNemirovski2015} working with stochastic subgradients into our framework is of great practical and theoretical interest. A critical assumption in our approach as well as others, e.g., see \cite{BenTalHazan2015} and references therein, is that the domain $X$ is convex.  Removing the convexity requirement on the domain $X$ will be an important theoretical development on its own. Besides, this will open up possibilities for more principled approaches to solving robust combinatorial optimization problems (see \cite{BertsimasBrownCaramanis2011,BertsimasSim2003}) where such a convexity assumption on $X$ is not satisfied. Finally, another attractive research direction is develop analogous frameworks for multi-stage RO problems  such as robust Markov decision processes (see \cite{NilimElGhaoui2005,Iyengar2005}).

\section*{Acknowledgments}
This research is supported in part by NSF grant CMMI 1454548.

The authors would like to thank Arkadi Nemirovski for suggesting the convex-concave reformulation of the robust feasibility problem presented in Section~\ref{sec:concaveSP}, in particular Lemma~\ref{lem:concaveSPReformulation}. The authors also would like to thank to the review team for useful suggestions and feedback that improved the presentation of the material in this paper.

\newpage
\appendix

\section{Convex-Concave Saddle Point Reformulation}\label{sec:concaveSP}

The SP problem~\eqref{eqn:convex-nonconcaveSP} based on the function $\Phi(x,u)$ which is not necessarily concave in $u$ admits a convex-concave SP representation in a lifted space via perspective transformations. To present this reformulation, we start by defining the following sets with additional variables $y\in\R^m_+$ and new variables $v^i$ for $i\in[m]$: 
\begin{align*}
V^i &= \left\{ [v^i;y^{(i)}] :~ 0 < y^{(i)} \leq 1,~ \frac{v^i}{y^{(i)}} \in U^i \right\}\quad\forall i\in[m], \\
W &= \left\{ w = [v^1;\ldots;v^m;y] :~ [v^i; y^{(i)}] \in \cl(V^i),\ i \in[m],~ \sum_{i=1}^{m} y^{(i)} = 1 \right\}.
\end{align*}
Note that for all $i\in[m]$, $\cl(V^i) = V^i \cup \{[0;0]\}$ because we assumed $U^i$ to be closed sets. 
For the point $[v^i;y^{(i)}] = [0;0]$, we set $y^{(i)} f^i\left( x, \frac{v^i}{y^{(i)}} \right) = 0$ for any $x\in X$. Note that setting $y^{(i)} f^i\left( x, \frac{v^i}{y^{(i)}} \right) = 0$ for  $[v^i;y^{(i)}] = [0;0]$ is well-defined as the continuation since from Assumption~\ref{ass:f-concave-u}, $f^i(x,u^i)$ is continuous and finite-valued on $U^i$, and $U^i$ is compact, so we deduce that $f^i(x,u^i)$ will be bounded on $U^i$. We also define the function $\psi:X \times W \to \R$ as
\[ 
\psi(x,w) = \psi(x,v,y) := \sum_{i=1}^{m} y^{(i)} f^i\left( x, \frac{v^i}{y^{(i)}} \right). 
\]

\begin{lemma}\label{lem:convexSPdef}
For fixed $w\in W$, the function $\psi(x,w)$ is convex in $x$ over $X$, and $\psi(x,w)$ is a concave function of $w$ over $W$ for any fixed $x$. Moreover, $W$ is closed, and when $U^i$ for $i\in[m]$ are convex, the sets $V^i$ for $i\in[m]$ and $W$ are all convex. 
\end{lemma}
\proof{Proof.}
For any $w = [v^1;\ldots;v^m;y]$, the function $\psi$ is convex in $x$ since in all of the nonzero terms in the summation over all $i\in[m]$ defining $\psi$, we have $y^{(i)}>0$ and in each such nonzero term each function $f^i\left( x, \frac{v^i}{y^{(i)}} \right)$ is convex in $x$ for the given $\frac{v^i}{y^{(i)}}\in U^i$ (see  Assumption~\ref{ass:f-concave-u}). In addition, for any given $x\in X$, the function $\psi$ is jointly concave in $v$ and $y$ because it is written as a sum of the perspective functions of functions $f^i$ which are concave in $u^i$ (see  Assumption~\ref{ass:f-concave-u}). 

The closedness of $W$ is immediate, and the convexity of the sets $V^i$ and $W$ follows immediately from their definition and the convexity assumption on $U^i$. 
\Halmos
\endproof

With these definitions and Lemma~\ref{lem:convexSPdef}, we observe that \eqref{eqn:convex-nonconcaveSP} is equivalent to evaluating the convex-concave SP problem defined by the function $\psi$ over the convex domains $X$ and $W$:
\begin{equation}\label{eqn:concaveSP}
\inf_{x \in X} \sup_{w \in W} \psi(x,w) \leq \epsilon \quad \text{or} \quad \inf_{x \in X} \sup_{w \in W} \psi(x,w) > 0.
\end{equation}
We state this formally in the following lemma. 
\begin{lemma}\label{lem:concaveSPReformulation}
For any $\epsilon > 0$ and $\bar{x} \in X$,
\[ 
\max_{i\in[m]} \sup_{u^i \in U^i} f^i(\bar{x},u^i) \leq \epsilon \quad\text{if and only if}\quad \sup_{w \in W} \psi(\bar{x},w) \leq \epsilon. 
\]
As a result,
\[ 
\inf_{x \in X} \max_{i\in[m]} \sup_{u^i \in U^i} f^i(x,u^i) \leq \epsilon \quad\text{if and only if}\quad \inf_{x \in X} \sup_{w \in W} \psi(x,w) \leq \epsilon. 
\]
\end{lemma}
\proof{Proof.} 
Fix $\bar{x} \in X$ and $\epsilon > 0$. Suppose $\max_{i\in[m]} \sup_{u^i \in U^i} f^i(\bar{x},u^i) \leq \epsilon$; then for all $u^i \in U^i$, $i \in [m]$, we have $f^i(\bar{x},u^i) \leq \epsilon$. Now consider any $w = [v^1;\ldots;v^m;y] \in W$. Then $0 \leq y^{(i)} \leq 1$ for all $i \in [m]$ and $\sum_{i=1}^m y^{(i)}=1$. For all $i \in [m]$, define $u^i = {v^i\over y^{(i)}} \in U^i$ whenever $y^{(i)}>0$.  
Then $y^{(i)} f^i(\bar{x},{v^i\over y^{(i)}}) = y^{(i)} f^i(\bar{x},u^i) \leq y^{(i)} \epsilon$ for $0 <y^{(i)} \leq 1$. In addition, when $y^{(i)}=0$, because $w\in W$ we must have $v^i=0$ and then by definition we have $y^{(i)} f^i(\bar{x},{v^i\over y^{(i)}})=0$. Therefore, from $\sum_{i=1}^my^{(i)}=1$, we deduce $\psi(\bar{x},w) = \sum_{i=1}^{m} y^{(i)} f^i\left(\bar{x}, \frac{v^i}{y^{(i)}} \right)\leq \epsilon$ holds for any $w\in W$.

Now suppose that $\sup_{w \in W} \psi(\bar{x},w) \leq \epsilon$ holds. Given $i \in [m]$ and $u^i \in U^i$, set $w$ to have components $y^{(i)} = 1$, $v^i = u^i$, and $[v^j;y^{(j)}] = [0;0]$ for $j \neq i$. Then $f^i(\bar{x},u^i) = \psi(\bar{x},w) \leq \epsilon$. Hence, $\max_{i\in[m]} \sup_{u^i \in U^i} f^i(\bar{x},u^i) \leq \epsilon$ follows.
\Halmos
\endproof

\begin{remark}\label{rem:m=1pers}
When $m=1$, i.e., we have only one function $f^1(x,u^1)$ and only one uncertainty set $U^1$, hence $W = U^1$ and $\inf_{x\in X} \sup_{v \in W} \psi(x,w) = \inf_{x\in X} \sup_{u^1 \in U^1} f^1(x,u^1)$. Also, under Assumption~\ref{ass:f-concave-u}, $\psi(x,w)$ is convex in $x$ and concave in $u^1$. Thus, the preceding perspective transformation resulting in \eqref{eqn:concaveSP} directly generalizes this case of a convex-concave SP formulation for $m=1$ discussed in Remark~\ref{rem:SPm=1}. 
\epr
\end{remark}

As a result, Lemma~\ref{lem:concaveSPReformulation} and Theorem~\ref{thm:saddlept-certificate} combined with any FOM that provides bounds on the saddle point gap $\epsilonsad^\psi(\bar{x},\bar{w})$ lead to an efficient way of verifying robust feasibility of \eqref{eqn:convex-nonconcaveSP} as follows: 
\begin{theorem}\label{thm:single-robust-constraint-feas-check}
Suppose $\bar{x}\in X$, $\bar{w}\in W$, and $\tau\in(0,1)$ are such that $\epsilonsad^\psi(\bar{x},\bar{w}) \leq \tau\epsilon$. If $\psi(\bar{x},\bar{w}) \leq (1-\tau)\epsilon$, then \sloppy $\max_{i \in [m]} \sup_{u^i \in U^i} f^i(\bar{x},u^i) \leq \epsilon$. If $\psi(\bar{x},\bar{w}) > (1-\tau)\epsilon$ and $\tau\leq{1\over 2}$, then $\inf_{x \in X} \max_{i \in [m]} \sup_{u^i \in U^i} f^i(x,u^i) > 0$.
\end{theorem}
\proof{Proof.}
Suppose $\psi(\bar{x},\bar{w}) \leq \tau\epsilon$. By Theorem~\ref{thm:saddlept-certificate}, we have $\inf_{x \in X} \sup_{w \in W} \psi(x,w) \leq \sup_{w\in W} \psi(\bar{x},w) \leq \epsilon$. By Lemma~\ref{lem:concaveSPReformulation}, $\max_{i \in [m]} \sup_{u^i \in U^i} f^i(\bar{x},u^i) \leq \epsilon$ as well.

On the other hand, when $\psi(\bar{x},\bar{w}) > (1-\tau)\epsilon$ and $\tau\leq{1\over 2}$, Theorem~\ref{thm:saddlept-certificate} implies \sloppy $\inf_{x\in X} \sup_{w\in W} \psi(x,w) \geq \inf_{x\in X} \psi(x,\bar{w}) > 0$. Then by Lemma~\ref{lem:concaveSPReformulation}, $\max_{i \in [m]} \sup_{u^i \in U^i} f^i(\bar{x},u^i) > 0$ follows.
\Halmos
\endproof

Because of the existence of efficient FOMs to solve convex-concave SP problems, Theorem~\ref{thm:single-robust-constraint-feas-check} suggests a possible advantage of using the convex-concave SP problem given in \eqref{eqn:concaveSP}. 
Nevertheless, working with the SP reformulation given by \eqref{eqn:concaveSP} in the extended space $X\times W$ presents a number of critical challenges. 
First, efficient FOMs associated with convex-concave SP problems often require computing prox operations or projections onto the domains $X$ and $W$. Unfortunately, even if projection (or prox-mappings) onto $U^i$ admits a closed form solution or an efficient procedure, it is unclear how to extend such projections onto $W$.
Furthermore, while the perspective transformations involved in constructing the function $\psi$ preserves certain desirable properties of the functions $f^i$, such as Lipschitz continuity and smoothness, the parameters associated with $\psi$ are in general larger than those associated with  the original functions $f^i$. Such parameters are critical for FOM convergence rates, and thus the FOMs when applied to solve \eqref{eqn:concaveSP} will have slower convergence rates. 

To address the issues outlined above, in the main paper we discuss how to obtain robust feasibility/infeasibility certificates for the convex-nonconcave SP problem \eqref{eqn:convex-nonconcaveSP} directly, i.e., we work with the functions $f^i$ and the sets $U^i$ directly. This direct approach in particular allows us to take greater advantage of the structure of the original formulation such as the availability of efficient projection (prox) computations over domains, and/or better parameters for smoothness, Lipschitz continuity, etc., of the functions.

\section{Supplementary Numerical Results}\label{sec:supplementary-numerical}
\renewcommand{\thefootnote}{\fnsymbol{footnote}}
Below, we provide the exact numerical values corresponding to the data used to generate the figures in our numerical study in Section \ref{sec:numerical}.

\begin{table}[ht]
\TABLE
{Average solve time (seconds) of each approach for different $n$ (Figure \ref{fig:average-solve-times}).
\label{tab:method-time-summary}}
{
\begin{tabular}{ll|cccc}
\toprule
    &  & \multicolumn{4}{c}{Approach} \\ \cline{3-6}
    &  & OFO-based & FO-based pessimization & nominal & full pessimization \\ \hline
\multirow{7}{*}{$n$} & 100 & 6.24\footnotemark[1] &  8.03\footnotemark[2] &     1.18 &      0.47 \\
    & 200 &         8.28\footnotemark[3] &           9.21\footnotemark[3] &     4.14 &      1.88 \\
    & 300 &        10.67 &          10.62 &     8.57 &      4.13 \\
    & 400 &        14.38\footnotemark[3] &          12.91\footnotemark[3] &    14.19 &      7.55 \\
    & 500 &        16.80 &          14.63 &    22.27 &     12.10 \\
    & 600 &        21.33 &          17.96 &    30.83 &     17.17 \\
    & 700 &        24.52 &          20.14 &    39.76 &     23.72 \\
\bottomrule
\end{tabular}
}
{}
\end{table}
\footnotetext[1]{Three instances out of 35 did not solve due to numerical issues.}
\footnotetext[2]{Two instances out of 35 did not solve due to numerical issues.}
\footnotetext[3]{One instance out of 35 did not solve due to numerical issues.}
\footnotetext[4]{One instance out of five did not solve due to numerical issues.}

\vspace{-5pt}
\begin{table}[ht]
\TABLE
{Average solve time (seconds) of each approach for different $m$ and $n$ (Figure \ref{fig:average-solve-times-plots}).
\label{tab:method-factors-time-summary}
}
{
\begin{tabular}{rll|rrrrrrr}
\toprule
    &             &   & \multicolumn{7}{c}{$m$} \\ \cline{4-10}
    &             &   &    3  &    5  &    7  &    10 &    15 &    20 &    25 \\ \hline
\multirow{16}{*}{$n$} & \multirow{4}{*}{400} & OFO-based & 13.17\footnotemark[4] & 12.21 & 14.78 & 14.03 & 14.73 & 16.23 & 15.24 \\
    &     & FO-based pessimization & 12.47\footnotemark[4] & 11.04 & 11.92 & 12.97 & 14.27 & 14.81 & 12.77 \\
    &     & nominal &  2.16 &  4.23 &  6.77 & 11.31 & 15.82 & 24.86 & 31.77 \\
    &     & full pessimization &  1.61 &  2.69 &  3.86 &  6.03 &  8.87 & 12.53 & 16.05 \\
\cline{2-10}
    & \multirow{4}{*}{500} & OFO-based & 15.96 & 15.55 & 15.88 & 17.22 & 17.41 & 17.52 & 18.04 \\
    &     & FO-based pessimization & 13.04 & 12.26 & 14.08 & 17.09 & 16.28 & 14.55 & 15.07 \\
    &     & nominal &  4.02 &  5.92 & 10.22 & 17.86 & 22.06 & 45.16 & 50.63 \\
    &     & full pessimization &  2.57 &  4.29 &  6.06 &  9.95 & 13.92 & 22.44 & 25.46 \\
\cline{2-10}
    & \multirow{4}{*}{600} & OFO-based & 18.62 & 19.52 & 20.97 & 21.25 & 20.86 & 23.46 & 24.62 \\
    &     & FO-based pessimization & 15.75 & 16.48 & 18.11 & 16.55 & 17.16 & 20.89 & 20.80 \\
    &     & nominal &  3.59 &  8.68 & 13.01 & 23.71 & 38.05 & 55.65 & 73.11 \\
    &     & full pessimization &  2.91 &  6.16 &  8.68 & 13.25 & 21.25 & 30.91 & 37.05 \\
\cline{2-10}
    & \multirow{4}{*}{700} & first-order & 19.49 & 21.86 & 24.46 & 25.04 & 25.12 & 28.48 & 27.18 \\
    &     & pessimization & 16.28 & 17.26 & 21.51 & 22.07 & 20.86 & 22.10 & 20.87 \\
    &     & nominal &  6.38 & 12.53 & 14.24 & 29.81 & 48.79 & 77.03 & 89.51 \\
    &     & full pessimization &  4.98 &  8.39 & 12.02 & 18.66 & 29.02 & 42.89 & 50.10 \\
\bottomrule
\end{tabular}
}
{}
\end{table}

\end{document}